\documentclass[11pt]{amsart}

\usepackage{amsmath,amsthm,amssymb}
\usepackage{hhline}
\usepackage{longtable}
\usepackage{mathrsfs}

\newtheorem{Le}{Lemma}[section]
\newtheorem{Pro}[Le]{Proposition}
\newtheorem{Th}[Le]{Theorem}
\newtheorem{Co}[Le]{Corollary}

\theoremstyle{definition}
\newtheorem{Def}[Le]{Definition}
\newtheorem{DaL}[Le]{Definition and Lemma}

\theoremstyle{remark}
\newtheorem{Ex}[Le]{Example}
\newtheorem{Rem}[Le]{Remark}

\allowdisplaybreaks
\begin{document}

\title[$G_2$-manifolds from K3 surfaces with a $\mathbb{Z}^2_2$-action]
{A construction of $G_2$-manifolds from K3 surfaces with a
$\mathbb{Z}^2_2$-action}
\author{Frank Reidegeld}
\address{TU Dortmund University,
Faculty for Mathematics,
44221 Dortmund, 
Germany}
\email{frank.reidegeld@math.tu-dortmund.de}

\begin{abstract}
A product of a K3 surface $S$ and a flat 3-dimensional torus $T^3$ 
is a manifold with holonomy $SU(2)$. Since $SU(2)$ is a subgroup
of $G_2$, $S\times T^3$ carries a torsion-free $G_2$-structure. We
assume that $S$ admits an action of $\mathbb{Z}^2_2$ with certain
properties. There are several possibilities to extend this action 
to $S\times T^3$. A recent result of Joyce and Karigiannis allows us 
to resolve the singularities of $(S\times T^3)/\mathbb{Z}^2_2$ such 
that we obtain smooth $G_2$-manifolds. We classify the quotients
$(S\times T^3)/\mathbb{Z}^2_2$ under certain restrictions and compute
the Betti numbers of the corresponding $G_2$-manifolds. Moreover,
we study a class of quotients by a non-abelian group. Several of our 
examples have new values of $(b^2,b^3)$. 
\end{abstract}

\maketitle

\tableofcontents

\section{Introduction}

A $G_2$-structure is a 3-form $\phi$ on a 7-dimensional manifold 
$M$ whose stabilizer group is $G_2$ at each point. Any 
$G_2$-structure induces a Riemannian metric. $(M,\phi)$ is called
a $G_2$-manifold if the holonomy group of the metric is $G_2$.
Currently, $G_2$-manifolds are an active area of research in pure 
mathematics and in theoretical physics \cite{Acharya}. Constructing 
new examples of compact $G_2$-manifolds turns out to be surprisingly
difficult. One of the reasons is that the equation for the holonomy 
reduction is highly non-linear. Another reason is that $G_2$-manifolds 
carry no complex structure and thus we cannot use techniques from
complex geometry as in the case of Calabi-Yau manifolds, which have 
holonomy $SU(n)$. 

There are only three known methods for the construction of compact 
$G_2$-manifolds. The first one, which led to the first compact examples 
of $G_2$-manifolds, is due to Dominic Joyce \cite{Joyce1,Joyce}. 
We choose a discrete group $\Gamma$ that acts on a torus $T^7$ and 
preserves the flat $G_2$-structure. It is possible to define a closed 
$G_2$-structure $\phi$ with small $\|d\ast\phi\|$ on a resolution of 
$T^7/\Gamma$. An analytic argument shows that $\phi$ can be 
deformed such that $d\ast\phi=0$ and the holonomy of the metric
thus is $G_2$ if the resolution is simply connected. 

The second construction is the twisted connected sum method, which
has been proposed by Donaldson and was carried out in detail by
Kovalev \cite{Ko}. Its starting point are two 
asymptotically cylindrical Calabi-Yau threefolds $W_1$ and $W_2$. 
Theorems that enable the construction of the $W_i$ can be found in 
\cite{HaskinsEtAl,Ko}. A $G_2$-manifold can be obtained by cutting off 
the cylindrical ends of the $W_i\times S^1$ and glueing together the truncated 
manifolds by an appropriate map. The first examples of 
$G_2$-manifolds that can be produced by this method can already 
be found in \cite{Ko}. In the meantime, several other classes of
examples have been constructed \cite{Braun,CortiEtAl1,CortiEtAl2,KoLe}. 

Recently, a third method was developed by Dominic Joyce and
Spiro Karigiannis \cite{JoyceKarigiannis}. Their idea is to divide a manifold
$M$ with a torsion-free $G_2$-structure $\phi$ by an involution
$\imath$. The quotient $M/\langle\imath\rangle$ has a singularity 
along a 3-dimensional submanifold $L$ that can be locally identified 
with $\mathbb{R}^3\times \mathbb{C}^2/\{\pm\text{Id}\}$. 
The next step is to cut out a neighbourhood of $L$ and glue in  
a family of Eguchi-Hanson spaces that are parameterized by $L$.   
After that, the authors obtain a smooth manifold that carries a metric 
whose holonomy is contained in $G_2$. It should be noted that the 
analytical methods that are needed to show the existence of a 
torsion-free $G_2$-structure are harder than in the two first 
constructions. One reason for this is that there is no canonical 
closed $G_2$-structure with small torsion on the family of 
Eguchi-Hanson spaces. Moreover, the method only works if a 
non-trivial condition is satisfied, namely that $L$ admits a closed 
and coclosed 1-form without zeroes. The result of \cite{JoyceKarigiannis}
can be easily generalized to the case where $M/\langle\imath\rangle$
is replaced by an arbitrary $G_2$-orbifold with singularities of
type $\mathbb{R}^3\times \mathbb{C}^2/\{\pm\text{Id}\}$.

The first two methods start with a manifold whose holonomy group
is smaller than $G_2$. The constructed $G_2$-manifold is in a certain 
sense nearby the manifold with smaller holonomy. In the first case, 
the starting point is a flat torus, which has trivial holonomy. The starting 
point of the twisted connected sum construction are the $W_i\times S^1$,
which have holonomy $SU(3)$. 

Beside the trivial group and $SU(3)$ there is also the group $SU(2)
\subseteq G_2$, which is also a possible holonomy group. Compact
$7$-dimensional manifolds with holonomy $SU(2)$ are covered by
$S\times T^3$, where $S$ is a K3 surface. This observation suggests
the following construction of $G_2$-manifolds. Since the holonomy of
$S\times T^3$ is a subgroup of $G_2$, it carries a torsion-free
$G_2$-structure $\phi$. Let $\Gamma$ be a discrete group that
acts on $S\times T^3$ and preserves $\phi$. The quotient
$(S\times T^3)/\Gamma$ is an orbifold. By resolving the singularities,
we should obtain a $G_2$-manifold. 

This idea has been already proposed by Joyce \cite{Joyce}. Moreover, 
he has shown that some of his resolutions of torus quotients are 
in fact examples of this method since K3 surfaces can be constructed 
as resolutions of $T^4/\mathbb{Z}_2$. As in \cite{JoyceKarigiannis},
the main obstacle to apply this method in full generality is to define a 
closed $G_2$-structure with small torsion on the resolved manifold. 
In the case $\Gamma=\mathbb{Z}_2$ we are in the situation of 
\cite{JoyceKarigiannis} and can work with their main result. 
Unfortunately, $(S\times T^3)/\mathbb{Z}_2$ has infinite 
fundamental group and therefore the holonomy of the manifold that we obtain
is not the whole group $G_2$. By dividing $S\times T^3$ by a group
that is isomorphic to $\mathbb{Z}^2_2$ and carefully applying
the theorem of \cite{JoyceKarigiannis} twice it is finally possible to obtain 
$G_2$-manifolds. 

Some examples of $G_2$-manifolds that are resolutions of 
$(S\times T^3)/\mathbb{Z}^2_2$ can be found in \cite{Joyce,JoyceKarigiannis}. 
The aims of this paper are to find larger classes of examples and to check if 
$G_2$-manifolds with the same Betti numbers can be found in the literature. 
We built up our research on the results of \cite{ReidK3} on K3 surfaces 
with a $\mathbb{Z}^2_2$-action. It turns out that we have to refine our 
results in order to obtain classifications of our examples. 

We consider three different extensions of the $\mathbb{Z}_2^2$-action
on $S$ to $S\times T^3$. Therefore, we obtain three different classes of 
examples. Each of them has distinct geometric properties. In the first case 
the action of one element of $\mathbb{Z}^2_2$ on $S\times T^3$ is free. 
This causes the holonomy of the manifold to be $SU(3)\rtimes\mathbb{Z}_2$ 
instead of $G_2$. In the second case our construction is equivalent to that 
of Kovalev and Lee \cite{KoLe}. Nevertheless, we find additional examples 
of $G_2$-manifolds that cannot be found in \cite{KoLe}. In the third case the 
second Betti number vanishes and $b^3$ takes only few different values. Finally, 
we consider another class of examples where $\Gamma$ is a non-abelian group 
instead of $\mathbb{Z}^2_2$. In that case we have $b^2,b^3\neq 0$ and obtain 
several $G_2$-manifolds that cannot be found in the existing literature.   

This paper is organised as follows. In Section \ref{G2-Section} we 
introduce the basic facts about $G_2$-manifolds and sum up the 
main theorems of \cite{JoyceKarigiannis}. In Section \ref{K3-Section} we recall the 
results of \cite{ReidK3} on K3 surfaces and refine them. Our examples of 
$G_2$-manifolds are contained in Section \ref{K3T3-Section}.

\section{$G_2$-manifolds and their construction methods}
\label{G2-Section}

\subsection{General facts} In this section, we introduce the most important facts
about $G_2$-manifolds that will be needed later on. For a more thorough 
introduction, we refer the reader to \cite{Joyce}. 

\begin{Def} Let $x^1,\ldots,x^7$ be the standard coordinates of $\mathbb{R}^7$
and let $dx^{i_1\ldots i_k} := dx^{i_1}\wedge\ldots\wedge dx^{i_k}$. The 3-form

\begin{equation}
\label{Can3form}
\phi_0 := dx^{123} + dx^{145} + dx^{167} + dx^{246} - dx^{257} - dx^{347} - dx^{356}
\end{equation}

is called the \emph{standard $G_2$-form on $\mathbb{R}^7$}.
\end{Def}

We define the group $G_2$ as the stabilizer of $\phi_0$ and mention some 
algebraic facts that are helpful to understand the inclusions $SU(2) \subseteq 
SU(3) \subseteq G_2$.

\begin{DaL}
\begin{enumerate}
    \item The group of all linear maps $f:\mathbb{R}^7\rightarrow \mathbb{R}^7$ with
    $f^{\ast}\phi_0=\phi_0$ is a connected Lie group whose Lie algebra is 
    the compact real of $\mathfrak{g}_2^{\mathbb{C}}$. We denote this group
    by $G_2$.
    
    \item Let $v\in\mathbb{R}^7\setminus\{0\}$. The group $\{f\in G_2|f(v)
    =v \}$ is isomorphic to $SU(3)$.

    \item A $3$-dimensional oriented subspace $V\subseteq\mathbb{R}^7$ is called
    \emph{associative} if $\phi_0|_V=\text{vol}$, where $\text{vol}$ is the 
    positive volume form on $V$ with respect to the standard Euclidean metric.
    The group of all $f\in G_2$ that act as the identity on $V$ is isomorphic 
    to $SU(2)$. 
\end{enumerate}    
\end{DaL}

Next, we introduce the notion of a $G_2$-structure. 

\begin{Def}
A \emph{$G_2$-structure} on a 7-dimensional manifold $M$ is a 3-form $\phi$ such
that for all $p\in M$ there exists a bijective linear map $f:T_pM\rightarrow \mathbb{R}^7$
such that $\phi_p = f^{\ast}\phi_0$.
\end{Def}

Since $G_2\subseteq SO(7)$, any $7$-dimensional manifold with a 
$G_2$-structure $\phi$ carries a canonical metric $g$. Moreover, the 
following equation defines a volume form $\text{vol}$ on $M$:

\begin{equation}
(v\lrcorner\phi)\wedge(w\lrcorner\phi)\wedge\phi = 6 g(v,w)\: \text{vol}\:.
\end{equation}

The following proposition helps us to decide if a $G_2$-structure induces a metric with
holonomy $G_2$. 

\begin{Pro}
Let $(M,\phi)$ be a manifold with a $G_2$-structure and let $g$ be the metric that
is induced by $\phi$. The following statements are equivalent. 

\begin{enumerate}
    \item $\nabla^g \phi = 0$, where $\nabla^g$ is the Levi-Civita connection.
    
    \item $d\phi = d\ast\phi = 0$.
    
    \item $\text{Hol} \subseteq G_2$, where $\text{Hol}$ is the holonomy group of
    the Levi-Civita connection. 
\end{enumerate}

If any of the above statements is true, $(M,g)$ is Ricci-flat. Conversely, let 
$(M,g)$ be a $7$-dimensional Riemannian manifold whose holonomy is a 
subgroup of $G_2$. If this is the case, $M$ carries a $G_2$-structure $\phi$
with $d\phi = d\ast\phi = 0$ whose induced metric is $g$.  
\end{Pro} 

\begin{Def}
In the situation of the above proposition, $\phi$ is called a \emph{torsion-free 
$G_2$-structure}. If in addition $\text{Hol} = G_2$, $(M,\phi)$ is called a 
\emph{$G_2$-manifold}. 
\end{Def}

If the underlying manifold is compact, it is particularly easy to decide if the 
holonomy is all of $G_2$ or just a subgroup.

\begin{Le}
Let $M$ be a compact manifold with a torsion-free $G_2$-structure. The 
holonomy of the induced metric is all of $G_2$ if and only if $\pi_1(M)$ is finite.   
\end{Le}

Finally, we have to introduce a special kind of submanifolds 
in order to formulate the results of \cite{JoyceKarigiannis}.

\begin{Def}
Let $(M,\phi)$ be a $G_2$-manifold and let $L$ be a 3-dimensional oriented 
submanifold of $M$. $L$ carries a volume form $vol$ that is determined by the
restriction of the induced metric to $L$ and its orientation. We call $L$
an \emph{associative submanifold of $M$} if for all $p\in L$ we have 
$\phi|_{T_pL}=vol_p$. 
\end{Def}

Since $G_2$ is the automorphism group of the octonions $\mathbb{O}$, 
a $G_2$-structure yields an identification of each tangent space with the 
imaginary space of the octonions $\text{Im}(\mathbb{O})$. An associative 
submanifold is characterized by the condition that each of its tangent spaces
is associative or equivalently that it can be identified with $\varphi(\text{Im}(
\mathbb{H}))$ for a $\varphi\in G_2$. Associative submanifolds are a special 
case of calibrated submanifolds. For a detailed introduction to this subject 
we refer the reader to the papers by Harvey and Lawson \cite{HarveyLawson} 
and by McLean \cite{McLean}.

\subsection{The construction of Joyce and Karigiannis} 
\label{JoyceKarigiannisSection} The main result of \cite{JoyceKarigiannis} can be stated as follows.

\begin{Th}
\label{Thm-Joyce-Karigiannis} (Theorem 6.4. in \cite{JoyceKarigiannis}) Let $(M,\phi)$
be a compact manifold with a torsion-free $G_2$-structure and let 
$\imath:M\rightarrow M$ be an involution that preserves $\phi$. We assume 
that $\imath$ is not the identity and has at least one fixed point. In this 
situation, the fixed point set of $\imath$ is a compact, not necessarily connected 
associative submanifold $L$. We assume that there exists a closed and
coclosed 1-form $\lambda$ on $L$ without zeroes.

The quotient $M/\langle \imath \rangle$ is a $G_2$-orbifold with an 
$A_1$-singularity along $L$, which means that all fibers of the normal bundle of 
$L$ can be identified with $\mathbb{C}^2/\{\pm\text{Id}\}$. With help of $\lambda$, 
it is possible to construct a bundle $\sigma:P \rightarrow L$ whose fibers are 
Eguchi-Hanson spaces. By cutting out a neighbourhood of $L\subseteq
M/\langle \imath \rangle$ and glueing in $P$, we obtain a resolution $\pi: N \rightarrow
M/\langle \imath \rangle$, where $N$ is a smooth manifold. 
The 5-dimensional manifold $\pi^{-1}(L)$ is a bundle over
$L$ whose fibers are diffeomorphic to $S^2$.  

In this situation, $N$ carries a torsion-free $G_2$-structure. In particular, $N$ is a
$G_2$-manifold if $\pi_1(N)$ is finite.  
\end{Th}

\begin{Rem}
\begin{enumerate}
    \item For the proof of the theorem it is not important that the 
    $A_1$-singularities are obtained by dividing $M$ by an involution. In
    fact, we may replace $M/\langle \imath \rangle$ by an arbitrary 
    $G_2$-orbifold with $A_1$-singularities.  

    \item For reasons of brevity, we call $N$ the resolution of $M/\langle\imath\rangle$, 
    although the process that is described in the above theorem is more complicated 
    than the resolution of singularities in algebraic geometry.  

    \item An Eguchi-Hanson space is the blow-up of $\mathbb{C}^2/\{\pm Id\}$
    together with a hyper-K\"ahler metric that approaches the standard 
    Hermitian metric on $\mathbb{C}^2/\{\pm Id\}$ at infinity. Since 
    $\mathbb{C}^2$ can be identified with $\mathbb{H}$, we have a sphere of 
    complex structures on $\mathbb{C}^2/\{\pm Id\}$. This sphere can be 
    identified with $S^2\subseteq \text{Im}(\mathbb{H})$. We are free to blow 
    up the singularity with respect to     any of these complex structures. The 
    value of the 1-form $\lambda$ at $p\in L$ is an element of $T^{\ast}_p L
    \setminus \{0\}$. Since $L$ is associative, the normalization of $\lambda_p$ 
    to unit length can be thought of as an element of $S^2\subseteq \text{Im}(
    \mathbb{H})$, too. Thus it determines a complex structure on the fiber of the 
    normal bundle. Therefore, it is possible to use $\lambda$ to define the fibers 
    of $P$. The details of this construction are important for the proof of Theorem 
    \ref{Thm-Joyce-Karigiannis}, but they do not influence the topology of the 
    $G_2$-manifold. We refer the reader to \cite{JoyceKarigiannis} for further information. 

    \item Although it is easy to decide if a closed and coclosed one-form 
    $\lambda$ exists on $L$, it is in general not possible to decide if $\lambda$ 
    has zeroes without information on the metric on $L$. Fortunately, the 
    submanifolds $L$ that we encounter are Riemannian products of  $S^1$ and a 
    two-dimensional manifold. Thus we can choose $\lambda = d\theta$ where 
    $\theta$ is the angle that parameterizes $S^1$.   
\end{enumerate}
\end{Rem}

The Betti numbers and the fundamental group of $N$ can be computed 
fairly easily. 

\begin{Co}
\label{JoyceKarigiannisBetti} (cf. Section 6.5 in \cite{JoyceKarigiannis}) In the situation of the above
two theorems we have 

\[
b^k(N) = b^k(M/\langle\imath\rangle) + b^{k-2}(L)
\] 

for all $k\in\{0,\ldots,7\}$, where $b^{-2}(L)$ and $b^{-1}(L)$ are defined as $0$. 
Moreover, the fundamental groups of $N$ and $M/\langle \imath \rangle$ are
isomorphic. 
\end{Co}

If we replace $\lambda$ with $-\lambda$, the complex structure $I$ on 
$\mathbb{C}^2/\{\pm Id\}$ is replaced by $-I$. Since this does not change the 
blow-up map and the sign of $\lambda$ plays no role in the other parts of the proof 
of Theorem \ref{Thm-Joyce-Karigiannis}, it suffices that $\lambda$ is defined up to a sign. We 
make this notion more explicit. Let $\pi:Z\rightarrow L$ be a $\mathbb{Z}_2$-principal
bundle over $L$. We can consider the bundles $\bigwedge^k T^{\ast} L 
\otimes_{\mathbb{Z}_2} Z$ over $L$. The fibers of these bundles are isomorphic to
$\bigwedge^k T^{\ast} L$ and its sections are called \emph{$Z$-twisted $k$-forms}.
If we have a local trivialization $L|_U \cong U\times \mathbb{Z}_2$
over an open subset $U\subset L$, we can identify $Z$-twisted $k$-forms naturally
with ordinary $k$-forms. The operators $d$ and $d^{\ast}$ induce operators on
the bundles $\bigwedge^k T^{\ast} L \otimes_{\mathbb{Z}_2} Z$ and therefore
it makes sense to talk about a closed and coclosed $Z$-twisted 1-form on $L$.
Moreover, we can define the \emph{$Z$-twisted de Rham cohomology} as the
cohomology of the complex $(\Gamma(\bigwedge^k T^{\ast} L \otimes_{\mathbb{Z}_2}
Z)_{k\geq 0},d)$ and the \emph{$Z$-twisted Betti numbers $b^k(L,Z)$} as the dimension
of the cohomologies. With this notation, there is the following generalization of    
Theorem \ref{Thm-Joyce-Karigiannis}. 

\begin{Co}
\label{Thm-Joyce-Karigiannis-Twisted} Let $(M,\phi)$, $\imath$ and $L$ be as
in Theorem \ref{Thm-Joyce-Karigiannis}. We assume that there exists a closed,
coclosed, non-vanishing $Z$-twisted 1-form $\lambda$ on $L$ for a  
$\mathbb{Z}_2$-bundle $Z$ on $L$. In this situation, the singularities of 
$M/\langle \imath \rangle$ can be resolved such that the resolved manifold $N$
carries a torsion-free $G_2$-structure. Moreover, we have 

\[
b^k(N) = b^k(M/\langle\imath\rangle) + b^{k-2}(L,Z)
\] 

for all $k\in\{0,\ldots,7\}$ and the fundamental groups of $N$ and $M/\langle 
\imath \rangle$ are isomorphic. 
\end{Co}

\section{K3 surfaces with non-symplectic $\mathbb{Z}^2_2$-actions}
\label{K3-Section}

\subsection{K3 surfaces with non-symplectic involutions}

We will see that any element of the group $\mathbb{Z}^2_2$, that 
acts on $S\times T^3$, acts on $S$ as a non-symplectic involution 
which is holomorphic with respect to one of the complex structures 
$I$, $J$ and $K$. Therefore, this section contains a brief introduction to 
K3 surfaces with non-symplectic involutions. A slightly more detailed account 
can be found in \cite{KoLe} or \cite{ReidK3}. We refer the reader for an 
in-depth introduction to the theory of K3 surfaces to \cite[Chapter VIII]{BHPV} 
and for the details of the classification of non-symplectic involutions to the 
original papers of Nikulin \cite{Nikulin1,Nikulin2,Nikulin3}. 

A K3 surface is a compact, simply connected, complex surface 
with tri\-vial canonical bundle. The underlying real $4$-dimensional manifold 
of a K3 surface is of a fixed diffeomorphism type. Therefore, any K3 surface 
$S$ has the same topological invariants. The Hodge numbers of $S$ are 
determined by $h^{0,0}(S)= h^{2,0}(S)=1$, $h^{1,0}(S)=0$ and 
$h^{1,1}(S)=20$. $H^2(S,\mathbb{Z})$ together with the intersection form
is an even, unimodular lattice with signature $(3,19)$. Up to isometries, 
the only lattice with these properties is the \emph{K3 lattice}

\begin{equation*}
L := 3 H \oplus 2(-E_8) \:,
\end{equation*}

where $H$ is the hyperbolic plane lattice with the quadratic form 
$q(x,y):=2xy$ and $-E_8$ is the root lattice of $E_8$ together with 
the negative of the usual bilinear form. We often write

\begin{equation*}
L = H_1 \oplus H_2 \oplus H_3 \oplus (-E_8)_1 \oplus (-E_8)_2
\end{equation*}

in order to distinguish between the different summands. We choose
a basis $(w_1,\ldots,w_{22})$ of $L$ such that $(w_{2j+1},
w_{2j+2})$ with $j=0,1,2$ is a basis of $H_j$ and

\begin{equation*}
w_{2j+1}\cdot w_{2j+1} = w_{2j+2}\cdot w_{2j+2} = 0\:,
\quad
w_{2j+1}\cdot w_{2j+2} = 1\:. 
\end{equation*}

Moreover, $(w_{7+8j},\ldots,w_{14+8j})$ with $j=0,1$ shall be a basis of 
$(-E_8)_j$ such that the matrix representation of the bilinear form with 
respect to  $(w_{7+8j},\ldots,w_{14+8j})$ is the negative of the Cartan 
matrix of $E_8$, which means that $w_k\cdot w_k=-2$ for $k\in\{7,
\ldots,22\}$ and for $k,l\in\{7,\ldots,22\}$ with $k\neq l$ we have
$w_k\cdot w_l \in \{0,1\}$. We call $(w_1,\ldots,w_{22})$ the 
\emph{standard basis of $L$}. 

Any K3 surface $S$ admits a K\"ahler metric. Since $S$ has trivial canonical 
bundle, there exists a unique Ricci-flat K\"ahler metric in each K\"ahler class. 
The holonomy group $SU(2)$ is isomorphic to $Sp(1)$. Therefore, the 
Ricci-flat K\"ahler metrics are in fact hyper-K\"ahler. In principle, a hyper-K\"ahler 
structure on a K3 surface is determined by the cohomology classes 
$[\omega_1],[\omega_2],[\omega_3]\in H^2(S,\mathbb{R})$ of the three 
K\"ahler forms. In order to make this statement precise, we introduce the
following notions.

\begin{Def} 
\begin{enumerate}
    \item Let $S$ be a K3 surface. A lattice isometry $\phi:H^2(S,\mathbb{Z})
    \rightarrow L$ is called a \emph{marking of} $S$. The pair $(S,\phi)$
    is called a \emph{marked K3 surface}.
    
    \item A \emph{hyper-K\"ahler structure on a marked K3 surface} is a tuple\linebreak
    $(S,\phi,g,\omega_1,\omega_2,\omega_3)$, where $g$ is a hyper-K\"ahler 
    metric and $\omega_i$ with $i\in\{1,2,3\}$ are the K\"ahler
    forms with respect to the complex structures $I_i$ that satisfy
    $I_1I_2I_3=-1$. We assume that $S$ has the orientation that makes $\omega_1
    \wedge \omega_1$ positive. 
    
    \item Two tuples $(S_1,\phi_1,g_1,\omega^1_1,\omega^1_2,\omega^1_3)$
    and $(S_2,\phi_2,g_2,\omega^2_1,\omega^2_2,\omega^2_3)$ are
    \emph{isometric} if there exists a map $f:S_1\rightarrow S_2$ with $f^{\ast}g_2=g_1$,
    $f^{\ast}\omega^2_i = \omega^1_i$ for $i\in\{1,2,3\}$ and $\phi_1\circ f^{\ast} = \phi_2$.     
    The \emph{moduli space of marked K3 surfaces with a hyper-K\"ahler
    structure} is the class of all tuples $(S,\phi,g,\omega_1,\omega_2,\omega_3)$
    modulo isometries.  
\end{enumerate}
\end{Def}

We write $L_{\mathbb{R}}$ for $L\otimes\mathbb{R}$ and define

\begin{equation*}
\begin{aligned}
\Omega:=\{ & (x,y,z)\in L_{\mathbb{R}}^3 | 
x^2 = y^2 = z^2>0,\: x\cdot y = x\cdot z = y\cdot z =0, \\
& \!\! \not\exists\: d\in L \:\:\text{with}\:\: d^2 = -2 \:\:\: \text{and} \:\:\:
x\cdot d = y\cdot d = z\cdot d=0    
\}\:. \\
\end{aligned}
\end{equation*}
    
The \emph{hyper-K\"ahler period map} $p$ that is defined by 
    
\[
p(S,\phi,g,\omega_1,\omega_2,\omega_3) = (\phi([\omega_1]),
\phi([\omega_2]),\phi([\omega_3]))
\] 

is a diffeomorphism between the hyper-K\"ahler moduli space 
and $\Omega$. This fact can also be found in \cite[Sec. 12.K]{Besse}
and \cite[Sec. 7.3]{Joyce}. The reason that we exclude triples $(x,y,z)$
with $x\cdot d = y\cdot d = z\cdot d=0$ from $\Omega$ is that 
they correspond to K3 surfaces with singularities. The following lemma
on isometries between K3 surfaces, that we have proven in \cite{ReidK3}, will
be useful later on. 

\begin{Le}
\label{IsomLem}
Let $S_j$ with $j\in\{1,2\}$ be K3 surfaces together with hyper-K\"ahler metrics $g_j$ 
and K\"ahler forms $\omega^j_i$ with $i\in\{1,2,3\}$. Moreover,
let $V_j\subseteq H^2(S_j,\mathbb{R})$ be the subspace that is spanned by 
the $[\omega^j_i]$, $i=1,2,3$.  

\begin{enumerate}
    \item Let $f:S_1\rightarrow S_2$ be an isometry. The pull-back
    $f^{\ast}:H^2(S_2,\mathbb{Z}) \rightarrow H^2(S_1,\mathbb{Z})$ is
    a lattice isometry. Its $\mathbb{R}$-linear extension maps $V_2$
    to $V_1$.

    \item \label{Isom} Let $\psi:H^2(S_2,\mathbb{Z}) \rightarrow 
    H^2(S_1,\mathbb{Z})$ be a lattice isometry. We denote the maps
    $\psi\otimes\mathbb{K}$ with $\mathbb{K}\in\{\mathbb{R},\mathbb{C}\}$
    by $\psi$, too. We assume that $\psi(V_2)=V_1$. Then there exists an 
    isometry $f:S_1\rightarrow S_2$ such that $f^{\ast} = \psi$. 

    \item Let $f:S\rightarrow S$ be an isometry that acts as the identity
    on $H^2(S,\mathbb{Z})$. Then, $f$ itself is the identity map. As a 
    consequence, the isometry from (\ref{Isom}) is unique.  
\end{enumerate}
\end{Le}

A \emph{non-symplectic involution} of a K3 surface $S$ is a holomorphic 
involution $\rho:S\rightarrow S$ such that $\rho$ acts as $-1$ on $H^{2,0}(S)$. 
Any K3 surface with a non-symplectic involution admits a K\"ahler class that is 
invariant under $\rho$. Therefore, there exists a hyper-K\"ahler structure on $S$ 
such that   

\[
\rho^{\ast}\omega_1 = \omega_1\:, 
\qquad 
\rho^{\ast}\omega_2 = - \omega_2\:,
\qquad
\rho^{\ast}\omega_3 = -\omega_3\:. 
\]   

Let $(S,\phi)$ be a marked K3 surface and $\rho:S\rightarrow S$
be a non-symplectic involution. The \emph{fixed lattice of $\rho$}
is defined as 

\begin{equation*}
L^{\rho} := \{x\in L | (\phi\circ\rho^{\ast}\circ\phi^{-1})(x) = x \}\:. 
\end{equation*}

$L^{\rho}$ is a non-degenerate sublattice of $L$ that is primitively
embedded, which means that $L/L^{\rho}$ has no torsion. Its 
signature is $(1,r-1)$, where $r$ is the rank of $L^{\rho}$. A 
lattice with this kind of signature is called \emph{hyperbolic}. 
Moreover, $L^{\rho}$ is \emph{2-elementary} which means that
$L^{\rho\ast}/L^{\rho}$ is isomorphic to a group of type 
$\mathbb{Z}_2^a$. The number $a\in\mathbb{N}_0$
is a second invariant of $L^{\rho}$. We define a third invariant
$\delta$ by

\begin{equation*}
\delta :=
\begin{cases}
0 & \text{if $x^2\in\mathbb{Z}$ for all $x\in L^{\rho\ast}$} \\
1 & otherwise \\
\end{cases}
\end{equation*} 

These invariants can in fact be defined for any $2$-elementary
lattice. If we assume that the lattice is in addition even and 
hyperbolic, there is at most one lattice with invariants $(r,a,\delta)$.  
Nikulin \cite{Nikulin3} has shown that the deformation classes of K3 
surfaces with a non-symplectic involution can be classified in terms 
of the triples $(r,a,\delta)$. There exist $75$ triples that satisfy

\begin{equation*}
1\leq r\leq 20\:,\quad 0\leq a\leq 11\quad\text{and}\quad
r-a\geq 0\:.
\end{equation*}

A figure with a graphical representation of all possible values
of $(r,a,\delta)$ can be found in \cite{KoLe, Nikulin3}. For the
computation of the Betti numbers of our $G_2$-manifolds we
need the following theorem about the fixed locus of a 
non-symplectic involution.

\begin{Th}
\label{FixedLocusTheorem}
(cf. \cite{KoLe,Nikulin3})
Let $\rho:S\rightarrow S$ be a non-symplectic involution of a K3 surface
and let $(r,a,\delta)$ be the invariants of its fixed lattice. The fixed locus 
$S^{\rho}$ of $\rho$ is a disjoint union of complex curves.

\begin{enumerate}
    \item If $(r,a,\delta)=(10,10,0)$, $S^{\rho}$ is empty.
    
    \item If $(r,a,\delta)=(10,8,0)$, $S^{\rho}$ is the disjoint union
    of two elliptic curves.
    
    \item In the remaining cases, we have
    
    \begin{equation*}
    S^{\rho} = C_g \cup E_1 \cup \ldots \cup E_k\:,
    \end{equation*}  
    
    where $C_g$ is a curve of genus $g=\frac{1}{2}(22 - r - a)$ and
    the $E_i$ are $k=\frac{1}{2}(r-a)$ rational curves.   
\end{enumerate}
\end{Th}

\subsection{Classification results} 
\label{NonSymplecticZ22}
For our construction of $G_2$-manifolds
we need K3 surfaces $S$ with a pair of commuting involutive isometries  
$(\rho^1,\rho^2)$ that satisfy the following equations:

\begin{equation}
\label{KahlerRelations1}
\begin{aligned}
& {\rho^1}^{\ast}\omega_1 = \omega_1\:, 
\qquad 
{\rho^1}^{\ast}\omega_2 = - \omega_2\:,
\qquad
{\rho^1}^{\ast}\omega_3 = -\omega_3 \\
& {\rho^2}^{\ast}\omega_1 = - \omega_1\:, 
\qquad 
{\rho^2}^{\ast}\omega_2 = \omega_2\:,
\qquad
{\rho^2}^{\ast}\omega_3 = - \omega_3 
\end{aligned}
\end{equation}

$\rho^3=\rho^1\rho^2$ is a third involution that satisfies

\[
{\rho^3}^{\ast}\omega_1 = - \omega_1\:, 
\qquad 
{\rho^3}^{\ast}\omega_2 = -\omega_2\:,
\qquad
{\rho^3}^{\ast}\omega_3 = \omega_3 
\] 

The set $(\text{Id},\rho^1,\rho^2,\rho^3)$ is a group that is isomorphic
to $\mathbb{Z}^2_2$. $\rho^i$ with $i=1,2,3$ preserves the hyper-K\"ahler
metric and the symplectic form $\omega_i$. Therefore, it is a biholomorphic 
map with respect to the complex structure $I_i$. Moreover, it maps the 
holomorphic $(2,0)$-form with respect to $I_i$ to its negative.
All in all, $\rho^i$ is a non-symplectic involution that is holomorphic with
respect to $I_i$. For this reason, we call the triple $(S,\rho^1,\rho^2)$ 
\emph{a K3 surface with a non-symplectic $\mathbb{Z}^2_2$-action}. 
Later on, we extend the $\rho^i$ with help of maps $\alpha^i:T^3
\rightarrow T^3$ to involutions of $S\times T^3$. Since the pairs
$(\rho^1\times \alpha^1,\rho^2\times \alpha^2)$ and
$(\rho^2\times \alpha^1,\rho^1\times \alpha^2)$ yield different quotients
$(S\times T^3)/\mathbb{Z}^2_2$, we consider $(\rho^1,\rho^2)$ and
$(\rho^2,\rho^1)$ as different objects. We should also keep in mind
that if $(\rho^1,\rho^2)$ generates a non-symplectic $\mathbb{Z}^2_2$-action,
the same is true for all $(\rho^i,\rho^j)$ with $i,j\in\{1,2,3\}$ and $i\neq j$. 

In order to find a classification result for non-symplectic 
$\mathbb{Z}^2_2$-actions that can be proven easily, we restrict ourselves 
to a special kind of non-symplectic involutions. 

\begin{Def}
\label{SimpleInvolution}
Let $S$ be a K3 surface and let $\rho:S\rightarrow S$  be
a non-symplectic involution. We call $\rho$ a \emph{simple
non-symplectic involution} if there exists a marking 
$\phi:H^2(S,\mathbb{Z})\rightarrow L$ such that for all
$i\in\{1,\ldots,22\}$ there exists a $j\in\{1,\ldots,22\}$
with $\rho(w_i)= \pm w_j$, where $\rho$ is an abbreviation
for $\phi\circ\rho^{\ast}\circ\phi^{-1}$ and $(w_1,\ldots,w_{22})$
is the standard basis of $L$.    
\end{Def} 

In \cite{ReidK3} we classified non-symplectic 
$\mathbb{Z}^2_2$-actions on K3 surfaces in terms of the
invariants $(r_i,a_i,\delta_i)$ of $\rho^i$, where $i\in\{1,2\}$. 
We made the assumptions that $\rho^1$ and $\rho^2$ are
simple and that the markings from the Definition 
\ref{SimpleInvolution} are the same for $\rho^1$ and $\rho^2$. 
In Section \ref{The_third_case}, we will need information on the
fixed loci of $\rho^1$, $\rho^2$ and $\rho^3$ in order to 
construct $G_2$-manifolds and to compute their Betti numbers.
Since the invariants $(r_3,a_3,\delta_3)$ 
are not determined by $(r_1,a_1,\delta_1)$ and
$(r_2,a_2,\delta_2)$, we need a refined version of our results
from \cite{ReidK3}. Since we make much use of the methods that we
have developed in \cite{ReidK3}, we recommend the reader to 
take a look at that paper.  
 
In \cite{ReidK3} we have shown that any simple non-symplectic 
involution $\rho$ maps an $H_i$ to an $H_j$ and an $(-E_8)_i$ to an 
$(-E_8)_j$. Since we have $\rho(w_i)= \pm w_j$ and each restriction 
$\rho|_{H_i}:H_i\rightarrow H_j$ has to be a lattice isometry, we can 
conclude after some calculations that each $\rho|_{H_i}$ is given 
by one of the following four matrices:

\[
M_1 := 
\begin{pmatrix}
1 & 0 \\
0 & 1 \\
\end{pmatrix}
\:, \quad
M_2 := 
\begin{pmatrix}
-1 & 0 \\
0 & -1 \\
\end{pmatrix}
\:, \quad
M_3 := 
\begin{pmatrix}
0 & 1 \\
1 & 0 \\
\end{pmatrix}
\:, \quad
M_4 := 
\begin{pmatrix}
0 & -1 \\
-1 & 0 \\
\end{pmatrix}
\]

Although the two matrices $M_3$ and $M_4$ are conjugate by a matrix
in $GL(2,\mathbb{R})$, the eigenspace of $M_3$ to the eigenvalue $1$ is 
spanned by a positive vector while the eigenspace to the same eigenvalue 
of $M_4$ is spanned by a negative vector. Therefore, $M_3$ and $M_4$ 
should be considered as two distinct cases. 

We denote the restriction of $\rho$ to $3H$ by $\rho':3H\rightarrow 3H$. 
If $\rho'(H_i)=H_i$ for all $i\in\{1,2,3\}$, $\rho'$ is a block matrix with 
$3$ blocks $M_j$, $M_k$ and $M_l$ along the diagonal. We denote
this diagonal block matrix by $\rho'_{jkl}$. Since the fixed lattice of 
$\rho$ has to be hyperbolic, one of the blocks is $M_1$ or $M_3$, 
which fix a positive vector, and the other ones have to be 
$M_2$ or $M_4$, which fix no positive vectors. Up to a permutation 
of the indices, $\rho'$ is one of the following 

\[
\rho'_{122}\:,\qquad 
\rho'_{124}\:,\qquad 
\rho'_{144}\:,\qquad 
\rho'_{322}\:,\qquad 
\rho'_{324}\:,\qquad 
\rho'_{344}\:. 
\]

Let $L_1,\ldots,L_n$ be $2$-elementary lattices with invariants 
$(r_i,a_i,\delta_i)$. The invariants of the direct sum $L_1\oplus
\ldots \oplus L_n$ are 

\[
(r_1+\ldots+r_n,a_1+\ldots+a_n,\max\{\delta_1,\ldots,\delta_n\})
\] 

We compute the fixed lattices of the $M_j$
and their invariants. After that, we obtain the invariants 
$(r',a',\delta')$ of the fixed lattices of the $\rho'_{jkl}$: 

\begin{center}
\label{table_rho_jkl}
\begin{tabular}{l|l|l}
$\rho|_{3H}$ & Fixed lattice & $(r',a',\delta')$ \\

\hline

$\rho'_{122}$ & $H$ & (2,0,0) \\
$\rho'_{124}$ & $H\oplus\mathbf{1}(-2)$ & (3,1,1) \\
$\rho'_{144}$ & $H\oplus\mathbf{1}(-2)\oplus\mathbf{1}(-2)$ & (4,2,1) \\
$\rho'_{322}$ & $\mathbf{1}(2)$ & (1,1,1) \\
$\rho'_{324}$ & $\mathbf{1}(2)\oplus\mathbf{1}(-2)$ & (2,2,1) \\
$\rho'_{344}$ & $\mathbf{1}(2)\oplus\mathbf{1}(-2)\oplus\mathbf{1}(-2)$ & (3,3,1) \\
\end{tabular}
\end{center} 

In the above table, $\mathbf{1}$ denotes the 1-dimensional lattice that
is generated by a single element of length $1$ and $L(k)$ with $k\in\mathbb{Z}$
denotes the lattice whose bilinear form is $k$ times the bilinear form of $L$. 

If there exist $i,j\in\{1,2,3\}$ with $i\neq j$ and $\rho'(H_i)=H_j$, the situation is more 
complicated. Since $\rho'$ is an involution, we have $\rho'(H_j)=H_i$ and the third 
sublattice of $3H$ is fixed. Without loss of generality, we assume that $\rho'(H_1)=H_2$, 
$\rho'(H_2)=H_1$ and $\rho'(H_3)=H_3$. Moreover, it follows from ${\rho'}^2=\text{Id}$
that the matrix representations of $\rho|_{H_1}$ and $\rho|_{H_2}$ have to be 
the same. This means that the matrix representation of $\rho'$ is given by

\[
\hat{\rho}'_{km} :=
\begin{pmatrix}
0 & M_k & 0 \\
M_k & 0 & 0 \\
0 & 0 & M_m \\
\end{pmatrix}   
\] 

Our next step is to determine the fixed lattice of $\hat{\rho}'_{km}$. We write
the elements of $3H$ as $(v_1,v_2,v_3)^{\top}$ with $v_i\in H_i$. By a short 
calculation we see that the fixed lattice is given by 

\[
\left\{
\begin{pmatrix}
v \\ M_k v \\ w
\end{pmatrix}
\middle|
v\in H,\: M_m w = w
\right\}
\]

The length of an element of the fixed lattice is

\[
2v^2 + w^2 
\] 

The fixed lattice contains a positive vector, namely $(v,M_kv,0)^{\top}$
with $v^2>0$. Since the fixed lattice has to be be hyperbolic, all vectors 
$w$ that are fixed by $M_m$ have to be negative. This implies that
$m\in \{2,4\}$, while $k\in\{1\ldots,4\}$ may be arbitrary. All in all,
we have found $8$ additional possibilities for $\rho'$. If $m=2$, the
fixed lattice of $\hat{\rho}'_{km}$ is isometric to $H(2)$ and has 
invariants $(2,2,0)$. If $m=4$, the fixed lattice is isometric to 
$H(2) \oplus \mathbf{1}(-2)$ and has invariants $(3,3,1)$. In the
case $m=4$, $\hat{\rho}'_{km}$ is actually conjugate to 
$\rho'_{344}$, but in the case $m=2$ it is not conjugate to any
$\rho'_{jkl}$ since the invariants $(2,2,0)$ do not appear in the
table from page \pageref{table_rho_jkl}. 

We denote the restriction of $\rho$ to $2(-E_8)$ by $\rho'':2(-E_8)
\rightarrow 2(-E_8)$. In \cite{ReidK3} we have shown that the 
restrictions $\rho|_{(-E_8)_i}:(-E_8)_i \rightarrow (-E_8)_j$ have to
be plus or minus the identity. Therefore, $\rho''$ has to be one of the 
matrices below.

\begin{align*}
& \rho''_1 = \begin{pmatrix}
I_8 & 0 \\
0 & I_8 \\
\end{pmatrix} \quad
\rho''_2 = \begin{pmatrix}
I_8 & 0 \\
0 & -I_8 \\
\end{pmatrix}  \quad
\rho''_3 = \begin{pmatrix}
-I_8 & 0 \\
0 & I_8 \\
\end{pmatrix} \\ 
& \rho''_4 = \begin{pmatrix}
-I_8 & 0 \\
0 & -I_8 \\
\end{pmatrix} \quad 
\rho''_5 = \begin{pmatrix}
0 & I_8 \\
I_8  & 0 \\
\end{pmatrix} \quad 
\rho''_6 = \begin{pmatrix}
0 & -I_8 \\
-I_8  & 0 \\
\end{pmatrix} \\
\end{align*} 

$\rho''_2$ is conjugate to $\rho''_3$ by a lattice isometry 
and the same is true for $\rho''_5$ and $\rho''_6$. No other pair of 
matrices from the above list is conjugate to each other. 
The invariants $(r'',a'',\delta'')$ of the fixed lattice of $\rho''$ can be found in 
the table below. 

\begin{center}
\begin{tabular}{l|l|l}
$\rho''$ & Fixed lattice & $(r'',a'',\delta'')$ \\

\hline

$\rho''_1$ & $2(-E_8)$ & (16,0,0) \\
$\rho''_2$ & $-E_8$ & (8,0,0) \\
$\rho''_4$ & $\{0\}$ & (0,0,0) \\
$\rho''_5$ & $-E_8(2)$ & (8,8,0) \\
\end{tabular}
\end{center} 

By writing $\rho$ as $\rho'\oplus\rho''$ and computing its invariants, we obtain our 
result from \cite{ReidK3}, which states that a non-symplectic involution is simple 
if and only if $(r,a,\delta)$ is an element of a list of $28$ triples.

We prove a classification result for non-symplectic $\mathbb{Z}^2_2$-actions. 
For reasons of simplicity, we make the following assumptions.

\begin{enumerate}
    \item The generators $\rho^i$ with $i=1,2$ of the group $\mathbb{Z}^2_2$
    are both simple involutions.
    
    \item According to Definition \ref{SimpleInvolution}, there exist markings
    $\phi^i:H^2(S,\mathbb{Z})\rightarrow L$ such that $\phi^i\circ {\rho^i}^{\ast}
    \circ (\phi^i)^{-1}$ has the desired matrix representation. In this article,
    we always assume that $\phi^1=\phi^2$ although further examples of
    non-symplectic $\mathbb{Z}^2_2$-actions may exist. 
    
    \item Let $L_i\subseteq L$ be the fixed lattice of $\rho^i$. We restrict 
    ourselves to the case that $L_1\cap L_2 = \{0\}$. We have proven in 
    \cite{ReidK3} that under this assumption we can find a smooth K3 surface 
    with a non-symplectic $\mathbb{Z}^2_2$-action and fixed lattices 
    $L_1$ and $L_2$ by choosing the three K\"ahler classes sufficiently 
    generic. 
\end{enumerate}

The $\rho^i$ can be written as ${\rho^i}'\oplus{\rho^i}''$ with ${\rho^i}': 3H
\rightarrow 3H$ and ${\rho^i}'':2(-E_8)\rightarrow 2(-E_8)$. We start with the
easier part and classify the pairs $({\rho^1}'',{\rho^2}'')$. Since ${\rho^i}''
\in \{\rho''_1,\ldots,\rho''_6 \}$, we have to check for all $j_1,j_2\in
\{1,\ldots,6\}$ if $\rho''_{j_1}$ and $\rho''_{j_2}$ commute and if the
sublattice of all vectors that are invariant under both maps is trivial.
If we count only those pairs $(\rho''_{j_1},\rho''_{j_2}) $ that cannot
be obtained from each other by conjugating both maps
$\rho''_{j_i}$ by an isometry of $2(-E_8)$, the following pairs $(j_1,j_2)$
remain.

\[ 
(1,4)\:,\quad 
(2,3)\:,\quad
(2,4)\:,\quad
(4,1)\:,\quad
(4,2)\:,\quad
(4,4)\:,\quad
(4,5)\:,\quad
(5,4)\:,\quad 
(5,6)
\]

Our next step is to classify the pairs $({\rho^1}',{\rho^2}')$. We start 
with the case where ${\rho^1}'$ and ${\rho^2}'$ map each $H_k$
to itself. Since the four matrices $M_1,\ldots,M_4$ commute pairwise,
any choice of the ${\rho^i}'$ yields commuting involutions. 
The fixed lattice of a non-symplectic involution is hyperbolic. 
Therefore, each ${\rho^i}'$ has to fix exactly one positive vector. 
Since we can permute the three spaces $H_1$, $H_2$ and $H_3$, we 
can assume without loss of generality that ${\rho^1}'|_{H_1}$ and 
${\rho^2}'|_{H_2}$ preserve a positive vector. These two restrictions 
have to be $M_1$ or $M_3$. The restrictions to the other $H_k$ have 
to be $M_2$ or $M_4$. As before, we have to take care that there 
are no lattice vectors in $3H$ that are fixed by both ${\rho^i}'$. The 
set of all pairs $(mn)$ with the property that the intersection of the 
eigenspaces of $M_m$ and $M_n$ to the eigenvalue $1$ is trivial 
consists of:

\[
(1,2) \quad (2,1)\quad (2,2)\quad (2,3)\quad (2,4)\quad (3,2)\quad 
(3,4)\quad (4,2)\quad (4,3)
\]  

With help of the above list, we find the following possibilities for 
${\rho^1}'$ and  ${\rho^2}'$. A triple $(p,q,r)$ in the $i$th column
of the table means that ${\rho^i}'=\rho'_{pqr}$. 

\begin{center}
\begin{tabular}{ccccc}
\begin{tabular}{l|l}
${\rho^1}'$ & ${\rho^2}'$ \\
\hline \hline
(1,2,2) & (2,1,2) \\
(1,2,2) & (2,1,4) \\
(1,2,2) & (2,3,2) \\
\hline
(1,2,2) & (2,3,4) \\
(1,2,4) & (2,1,2) \\
(1,2,4) & (2,3,2) \\
\hline
(1,4,2) & (2,3,2) \\
(1,4,2) & (2,3,4) \\
(1,4,4) & (2,3,2) \\
\end{tabular}
& $\qquad$ & 
\begin{tabular}{l|l}
${\rho^1}'$ & ${\rho^2}'$ \\
\hline \hline
(3,2,2) & (2,1,2) \\
(3,2,2) & (2,1,4) \\
(3,2,2) & (2,3,2) \\
\hline
(3,2,2) & (2,3,4) \\
(3,2,2) & (4,1,2) \\
(3,2,2) & (4,1,4) \\
\hline
(3,2,2) & (4,3,2) \\
(3,2,2) & (4,3,4) \\
(3,2,4) & (2,1,2) \\
\end{tabular}
& $\qquad$ & 
\begin{tabular}{l|l}
${\rho^1}'$ & ${\rho^2}'$ \\
\hline \hline 
(3,2,4) & (2,3,2) \\
(3,2,4) & (4,1,2) \\
(3,2,4) & (4,3,2) \\
\hline
(3,4,2) & (2,3,2) \\
(3,4,2) & (2,3,4) \\
(3,4,2) & (4,3,2) \\
\hline
(3,4,2) & (4,3,4) \\
(3,4,4) & (2,3,2) \\
(3,4,4) & (4,3,2) \\
\end{tabular} \\
\end{tabular}
\end{center}

\begin{Rem}
\label{RemarkLatticeInvariants}
In the above table, there are a few cases where the invariants of 
${\rho^1}'$ and ${\rho^2}'$ are the same but the $\mathbb{Z}^2_2$-actions
cannot be obtained from each other by conjugating with a lattice isometry. 
An example can be found in the 6th and 7th row of the table. The fixed 
lattice of ${\rho^2}'$ is in both cases $\mathbf{1}(2)$ which is embedded 
into $H_2$. The fixed lattice of ${\rho^1}'$ is $H_1\oplus \mathbf{1}(-2)$, 
but in the 6th line $\mathbf{1}(-2)$ is embedded into $H_3$ and in the 7th 
line it is embedded into $H_2$. In the 6th line $\mathbf{1}(-2)$ is embedded
into a different $H_k$ as $\mathbf{1}(2)$ and in the 7th line they are both 
embedded into $H_2$. Therefore, the $\mathbb{Z}^2_2$-actions are 
inequivalent although the fixed lattices are isomorphic. Since it is not clear 
if we can obtain topologically different $G_2$-manifolds from those two 
pairs, we will consider them as two distinct cases. 
\end{Rem}

Next, we assume that ${\rho^1}'$ interchanges two of the $H_k$ and 
${\rho^2}'$ leaves all $H_k$ invariant. Without loss of generality 
${\rho^1}'$ interchanges $H_1$ and $H_2$. We recall that ${\rho^1}'$ 
is a matrix of type 

\[
\hat{\rho}'_{km} = 
\begin{pmatrix}
0 & M_k & 0 \\
M_k & 0 & 0 \\
0 & 0 & M_m \\
\end{pmatrix}   
\] 

with $k\in\{1,\ldots,4\}$ and $m\in\{2,4\}$. $\hat{\rho}'_{km}$ fixes exactly
one positive vector and this vector is an element of $H_1\oplus H_2$. By a short 
calculation we see that a diagonal block matrix $\rho'_{pqr}$ commutes 
with $\hat{\rho}'_{km}$ if and only if the following three equations are satisfied

\begin{eqnarray*}
M_k M_q & = & M_p M_k \\
M_k M_p & = & M_q M_k \\
M_m M_r & = & M_r M_m \\
\end{eqnarray*}  

Since the matrices $M_1,\ldots,M_4$ commute pairwise, the third equation 
is automatically satisfied. The first two equations are equivalent to  

\[
M_q = M_k M_p M_k 
\]

since $M_k^2=1$. Again we use the fact that $M_1,\ldots,M_4$ commute
and conclude that $p=q$. A linear map $\rho'_{ppr}$ fixes exactly one
positive vector if and only if $p\in\{2,4\}$ and $r\in\{1,3\}$. In all cases, the fixed 
vector is an element of $H_3$ and thus the intersection of the fixed lattices is
trivial. All in all, there are $8$ possibilities for ${\rho^1}'$ and $4$ possibilities for 
${\rho^2}'$. Therefore, we have found $32$ additional pairs $({\rho^1}',{\rho^2}')$. 
We conjugate both maps by the matrix

\[
\begin{pmatrix}
1 & 0 & 0 \\
0 & M_k & 0 \\
0 & 0 & 1 \\
\end{pmatrix}
\]

This conjugation transforms ${\rho^1}'$ into $\hat{\rho}'_{1m}$ and
leaves ${\rho^2}'$ invariant. Therefore, only the following $8$ pairs 
remain:

\[
\begin{tabular}{llll}
$(\hat{\rho}'_{12},\rho'_{221})\:,\quad$ &
$(\hat{\rho}'_{12},\rho'_{223})\:,\quad$ &
$(\hat{\rho}'_{12},\rho'_{441})\:,\quad$ &
$(\hat{\rho}'_{12},\rho'_{443})\:,$ \\[2mm]
$(\hat{\rho}'_{14},\rho'_{221})\:,\quad$ &
$(\hat{\rho}'_{14},\rho'_{223})\:,\quad$ &
$(\hat{\rho}'_{14},\rho'_{441})\:,\quad$ &
$(\hat{\rho}'_{14},\rho'_{443})$ \\
\end{tabular}
\]

In the case where ${\rho^1}'$ leaves all $H_k$ invariant and ${\rho^2}'$ 
interchanges two of these lattices, we find $8$ analogous pairs. 
Finally, we consider the case where both ${\rho^1}'$ 
and ${\rho^2}'$ interchange two summands of $3H$. Since 
${\rho^1}'$ and ${\rho^2}'$ shall commute, this is only possible
if both of them interchange the same pair. Therefore, we assume
without loss of generality that ${\rho^1}'$ and ${\rho^2}'$ interchange 
$H_1$ and $H_2$. We have ${\rho^1}' =  \hat{\rho}'_{km}$ and
${\rho^2}' =  \hat{\rho}'_{ln}$ with $k,l\in \{1,\ldots,4\}$ and $m,n\in \{2,4\}$.
It is easy to see that ${\rho^1}'$ and ${\rho^2}'$ commute, no matter
what the values of $k,l,m$ and $n$ are. Therefore, we have found further 
$64$ pairs $({\rho^1}',{\rho^2}')$. By conjugating with the same matrix as in 
the previous case we can achieve that $k=1$, but $l,m$ and $n$ can still be 
chosen arbitrarily. This reduces the number of pairs to $16$. 

We obtain the complete list of all simple non-symplectic $\mathbb{Z}^2_2$-actions
that satisfy our additional assumptions $\phi_1=\phi_2$ and $L_1\cap L_2 = \{0\}$
by combining any $({\rho^1}',{\rho^2}')$ with any $({\rho^1}'',{\rho^2}'')$. All in
all, we find $(27 + 8 + 8 + 16)\cdot 9 = 531$ pairs. 

We denote the invariants of $\rho^i$ with $i\in\{1,2,3\}$ by $(r_i,a_i,\delta_i)$. 
Later on, we will see that the Betti numbers of our $G_2$-manifolds are 
independent of the $\delta_i$. Since in one of our cases the Betti numbers 
depend on the invariants of all $\rho^i$, we search for all tuples 
$(r_1,a_1|r_2,a_2|r_3,a_3)$ that can be obtained by our construction. 
We have $(r_i,a_i)=(r'_i + r''_i,a'_i+a''_i)$, where $(r'_i,a'_i)$ are the invariants 
of ${\rho^i}'$ and $(r''_i,a''_i)$ are the invariants of ${\rho^i}''$. We go through 
our lists of pairs $({\rho^1}',{\rho^2}')$ and $({\rho^1}'',{\rho^2}'')$, determine 
${\rho^3}'$ and ${\rho^3}''$ and compute the invariants of the lattice involutions. 
After that, we have finally proven the following theorem. 

\begin{Th}
\label{Theorem_Z22} 
Let $(\rho^1,\rho^2)$ be a pair of commuting non-symplectic involutions
of a smooth K3 surface that are holomorphic with respect to complex structures
$I_1$ and $I_2$ with $I_1I_2=-I_2I_1$. In addition, we assume that

\begin{enumerate}
    \item the intersection of the fixed lattices of $\rho^1$ and $\rho^2$ is trivial,

    \item $\rho^1$ and $\rho^2$ are simple and
    
    \item the markings from Definition \ref{SimpleInvolution} are the same for
    $\rho^1$ and $\rho^2$. 
\end{enumerate}

We denote $\rho^1\rho^2$ by $\rho^3$. The action of each $\rho^i$ 
on the K3 lattice $L$ splits into ${\rho^i}'\oplus {\rho^i}''$ with  
${\rho^i}':3H\rightarrow 3H$ and ${\rho^i}'':2(-E_8) \rightarrow 2(-E_8)$. 
Moreover, the action of $\rho^1$ and $\rho^2$ is conjugate by a lattice
isometry to one of $531$ actions that we have described earlier in this section.   

We denote the invariants of $\rho^i$ by $(r_i,a_i,\delta_i)$. The set of all 
possible tuples $(r_1,a_1|r_2,a_2|r_3,a_3)$ can be obtained by
adding a tuple from the list of all invariants $(r'_1,a'_1|r'_2,a'_2|r'_3,a'_3)$
of the ${\rho^i}'$:

\begin{center}
\begin{tabular}{lll}
$(1,1|1,1|4,2)$ & $(1,1|4,2|1,1)$ & $(4,2|1,1|1,1)$ \\ [2mm]
$(1,1|2,0|3,1)$ & $(1,1|3,1|2,0)$ & $(2,0|1,1|3,1)$ \\ 
$(2,0|3,1|1,1)$ & $(3,1|1,1|2,0)$ & $(3,1|2,0|1,1)$ \\ [2mm]
$(1,1|2,2|3,1)$ & $(1,1|3,1|2,2)$ & $(2,2|1,1|3,1)$ \\
$(2,2|3,1|1,1)$ & $(3,1|1,1|2,2)$ & $(3,1|2,2|1,1)$ \\ [2mm]
$(1,1|2,2|3,3)$ & $(1,1|3,3|2,2)$ & $(2,2|1,1|3,3)$ \\
$(2,2|3,3|1,1)$ & $(3,3|1,1|2,2)$ & $(3,3|2,2|1,1)$ \\ [2mm]
$(2,0|2,0|2,0)$ & & \\ [2mm]
$(2,0|2,2|2,2)$ & $(2,2|2,0|2,2)$ & $(2,2|2,2|2,0)$ \\ [2mm]
$(2,0|3,3|3,3)$ & $(3,3|2,0|3,3)$ & $(3,3|3,3|2,0)$ \\ [2mm]
$(2,2|2,2|2,2)$ & & \\ [2mm]
$(2,2|2,2|4,2)$ & $(2,2|4,2|2,2)$ & $(4,2|2,2|2,2)$ \\ [2mm] 
$(2,2|3,3|3,3)$ & $(3,3|2,2|3,3)$ & $(3,3|3,3|2,2)$ \\ [2mm]
$(3,3|3,3|4,2)$ & $(3,3|4,2|3,3)$ & $(4,2|3,3|3,3)$ \\
\end{tabular}
\end{center}

to a tuple from the list of all invariants $(r''_1,a''_1|r''_2,a''_2|r''_3,a''_3)$
of the ${\rho^i}''$:

\begin{center}
\begin{tabular}{lll}
$(0,0|0,0|16,0)$ & $(0,0|16,0|0,0)$ & $(16,0|0,0|0,0)$ \\
$(0,0|8,0|8,0)$ & $(8,0|0,0|8,0)$ & $(8,0|8,0|0,0)$ \\
$(0,0|8,8|8,8)$ & $(8,8|0,0|8,8)$ & $(8,8|8,8|0,0)$ \\
\end{tabular}  
\end{center}

In particular, there are $38\cdot 9 = 342$ different tuples 
$(r_1,a_1|r_2,a_2|r_3,a_3)$. 
\end{Th}

Another large class of non-symplectic $\mathbb{Z}^2_2$-actions is
implicitly contained in the article of Kovalev and Lee \cite{KoLe} on 
twisted connected sums. Their construction of $G_2$-manifolds
can be summed up as follows.

\begin{enumerate}
    \item Let $S$ be a K3 surface and $\rho:S\rightarrow S$ be a 
    non-symplectic involution. Moreover, let $\psi$ be a holomorphic 
    involution of $\mathbb{CP}^1$. 

    \item The quotient $Z:=(S\times \mathbb{CP}^1)/\langle\rho\times
    \psi\rangle$ is an orbifold with $A_1$-singularities along two copies
    of the fixed locus of $\rho$ since $\psi$ has two fixed points. 

    \item We blow up the singularities of $Z$ and obtain a compact
    K\"ahler manifold $\overline{W}$. 

    \item $\overline{W}$ is fibered by K3 surfaces. We remove a 
    fiber from $\overline{W}$ and obtain a non-compact manifold
    $W$ that carries an asymptotically cylindrical Ricci-flat K\"ahler
    metric that approaches $S\times S^1\times (0,\infty)$. 

    \item Let $W_1$ and $W_2$ be two asymptotically cylindrical 
    Calabi-Yau manifolds that can be obtained by the above construction. 
    $W_1\times S^1$ and $W_2\times S^1$ can be glued together
    to a $G_2$-manifold if a \emph{matching} between $S_1$
    and $S_2$ exists, that is an isometry $f:S_1\rightarrow S_2$ 
    with 
    
    \[
    f^{\ast}[\omega_1^2] = [\omega_2^1]\:,\qquad
    f^{\ast}[\omega_2^2] = [\omega_1^1]\:,\qquad
    f^{\ast}[\omega_3^2] = -[\omega_3^1]\:.
    \]
    
    In the above equation, $[\omega_j^i]$ denotes the cohomology
    class of the $j$th K\"ahler form on the K3 surface $S_i$.
\end{enumerate}  

The first three steps describe the main idea of Kovalev and Lee
\cite{KoLe}. The fourth step is a general result on asymptotically 
cylindrical Calabi-Yau manifolds that can be found in 
\cite{HaskinsEtAl,Ko} and the fifth step is the actual twisted connected 
sum construction that is developed in \cite{Ko}. The authors show 
with help of a theorem on lattice embeddings that a matching exists 
if the invariants $(r_i,a_i,\delta_i)$ of the non-symplectic involutions 
$\rho^i$ of $S_i$ satisfy $r_1+r_2\leq 11$ or $r_1+r_2+a_1+a_2<22$. 
The above construction yields a large number of $G_2$-manifolds. 
Further details including the Betti numbers of the $G_2$-manifolds 
can be found in \cite{KoLe}. Since $\rho^1$ and $f^{-1} \circ \rho^2 
\circ f$ define a non-symplectic $\mathbb{Z}^2_2$-action on $S_1$, 
we immediately obtain the following theorem.

\begin{Th}
\label{KoLe-Theorem}
Let $(r_i,a_i,\delta_i)\in \mathbb{N}\times\mathbb{N}_0\times \{0,1\}$ 
with $i=1,2$ be triples such that non-symplectic involutions of
K3 surfaces with invariants $(r_i,a_i,\delta_i)$ exist. If

\[
r_1+r_2\leq 11 \quad\text{or}\quad r_1+r_2+a_1+a_2<22
\]

a K3 surface with a non-symplectic $\mathbb{Z}^2_2$-action 
$(S,\rho^1,\rho^2)$ exists such that $\rho^i$ has invariants 
$(r_i,a_i,\delta_i)$. 
\end{Th}

\section{Examples of $G_2$-manifolds}
\label{K3T3-Section}

\subsection{The idea behind our construction} 
Let $S$ be a K3 surface with a hyper-K\"ahler structure and let $T^3 = 
\mathbb{R}^3/\Lambda$, where $\Lambda$ is a lattice of rank $3$, be 
a flat torus. The metric on the product $S\times T^3$ has holonomy 
$Sp(1)\subseteq G_2$ and $S\times T^3$ therefore carries a torsion-free 
$G_2$-structure, which we describe in detail. Let $x^1$, $x^2$ and $x^3$ be 
coordinates on $T^3$ such that $(\tfrac{\partial}{\partial x^1},
\tfrac{\partial}{\partial x^2},\tfrac{\partial}{\partial x^3})$ is an orthonormal 
frame and let $\omega_1,\omega_2,\omega_3$ be the three K\"ahler forms 
on $S$. The 3-form

\begin{equation*}
\phi:= \omega_1\wedge dx^1 + \omega_2\wedge dx^2 + 
\omega_3\wedge dx^3 + dx^1\wedge dx^2\wedge dx^3 
\end{equation*}

is a torsion-free $G_2$-structure, whose associated metric is the product 
metric on $S\times T^3$. Let $\Gamma$ be a finite group that acts 
on $S\times T^3$ and leaves $\phi$ invariant. 
The action of any $\gamma\in\Gamma$ can be written as a 
product of an isometry of $S$ and an isometry of $T^3$. Since 
$\phi$ contains the summand $dx^1\wedge dx^2\wedge dx^3$, 
$\gamma$ has to preserve the orientation of $T^3$. The action 
of $\gamma$ on $T^3$ can thus be written as

\begin{equation*}
x + \Lambda \mapsto (A^{\gamma} x + v^{\gamma}) + \Lambda\:, 
\end{equation*}

where $v^{\gamma}\in\mathbb{R}^3$ and $A^{\gamma} = 
(A_{ij}^{\gamma})_{i,j=1,2,3}\in SO(3)$. Let $GL(\Lambda)$ 
be the subgroup of $GL(3,\mathbb{R})$ that maps 
$\Lambda$ to itself. $\gamma$ defines a well-defined, orientation 
preserving isometry of $T^3$ if and only if $A^{\gamma} \in 
SO(3)\cap GL(\Lambda)$. We assume that for any $\gamma
\in\Gamma$ there exists an isometry of $S$ whose pull-back
acts on the three K\"ahler forms as 

\begin{equation*}
\omega_i \mapsto \sum_{j=1}^3 A_{ji}^{\gamma}\omega_j\:. 
\end{equation*}

In this case the action of $\Gamma$ on $T^3$ can be extended 
to an action on $S\times T^3$ that leaves $\phi$ invariant. If 
$\Gamma$ acts freely on $S\times T^3$, the quotient
$(S\times T^3)/\Gamma$ is a manifold whose holonomy group has 
$Sp(1)$ as identity component. In the following subsection, we 
search for non-free group actions on $S\times T^3$ such that 
we can resolve the singularities of $(S\times T^3)/\Gamma$ by 
Theorem \ref{Thm-Joyce-Karigiannis} and obtain smooth 
$G_2$-manifolds.

\subsection{The choice of the group action} In order to motivate how
the group $\Gamma$ acts on $S\times T^3$, we recall what is 
known about free group actions on $T^3$. A quotient of a torus by an 
isometric free group action is called a compact Euclidean space form. 
In dimension $3$, they were classified by Hantzsche and Wendt 
\cite{Hantzsche}. In the case where $\Gamma$ preserves the 
orientation there exist $6$ space forms. Only the last one, which 
has finite fundamental group, is relevant for our construction. 
In that case, the lattice $\Lambda$ can be chosen as $\mathbb{Z}^3$ 
and $\Gamma$ is isomorphic to $\mathbb{Z}_2^2$. The generators
$\psi^1$ and $\psi^2$ of $\mathbb{Z}_2^2$ are given by 

\begin{equation}
\begin{array}{rcl}
\label{Generators}
\psi^1((x^1,x^2,x^3) + \mathbb{Z}^3) & := & 
(\tfrac{1}{2}+x^1, -x^2,\tfrac{1}{2} - x^3) + \mathbb{Z}^3 \\
\psi^2((x^1,x^2,x^3) + \mathbb{Z}^3) & := & 
(-x^1,\tfrac{1}{2} + x^2, -x^3) + \mathbb{Z}^3 \\
\end{array}
\end{equation}

We extend the maps $\psi^i$ with $i=1,2$ to maps $\rho^i\times \psi^i:
S\times T^3 \rightarrow S\times T^3$ that leave $\phi$ invariant and
generate a $\mathbb{Z}^2_2$-action. It follows from (\ref{Generators}) 
that the pull-backs of $\rho^1$ and $\rho^2$ have to act on the K\"ahler 
forms as in the equation (\ref{KahlerRelations1}) that we have studied in 
Section \ref{NonSymplecticZ22}. 

$\mathbb{Z}_2^2$ acts freely on $T^3$ and 
$(S\times T^3)/\mathbb{Z}_2^2$ thus is a smooth manifold with holonomy 
$Sp(1) \rtimes \mathbb{Z}_2^2$. Since this is not what we want, we modify the 
translation part of $\psi^1,\psi^2:T^3\rightarrow T^3$ such that the action of
$\mathbb{Z}_2^2$ on $T^3$ is not free anymore, but the action of 
$\mathbb{Z}_2^2$ on $S$ still satisfies the relations (\ref{KahlerRelations1}). 
All elements of $\mathbb{Z}_2^2$ are of order $2$ and we will see that the 
singularities of $(S\times T^3)/\mathbb{Z}_2^2$ look locally like 
$\mathbb{C}^2/\{\pm 1\} \times \mathbb{R}^3$ and thus can be resolved by 
Theorem \ref{Thm-Joyce-Karigiannis}. We consider three different kinds of 
group actions. In the first case we replace the translation part of $\psi^2$ 
by zero and obtain: 

\begin{equation}
\begin{array}{rcl}
\label{K3T3Case1}
\psi^1((x^1,x^2,x^3) + \mathbb{Z}^3) & := & (\tfrac{1}{2}+x^1,-x^2,\tfrac{1}{2} - x^3) + \mathbb{Z}^3 \\
\psi^2((x^1,x^2,x^3) + \mathbb{Z}^3) & := & (-x^1,x^2,-x^3) + \mathbb{Z}^3 \\
\psi^3((x^1,x^2,x^3) + \mathbb{Z}^3) & := & (\tfrac{1}{2}-x^1,-x^2,\tfrac{1}{2}+x^3) + \mathbb{Z}^3 \\
\end{array}
\end{equation}

where $\psi^3:=\psi^1\psi^2$. In the second case we do the same 
for $\psi^1$ and have: 

\begin{equation}
\begin{array}{rcl}
\psi^1((x^1,x^2,x^3) +  \mathbb{Z}^3) & := & (x^1,-x^2,- x^3) +  \mathbb{Z}^3 \\
\psi^2((x^1,x^2,x^3) +  \mathbb{Z}^3) & := & (-x^1,\tfrac{1}{2}+x^2,-x^3) +  \mathbb{Z}^3 \\
\psi^3((x^1,x^2,x^3) +  \mathbb{Z}^3) & := & (-x^1,\tfrac{1}{2}-x^2,x^3) +  \mathbb{Z}^3 \\
\end{array}
\end{equation}

For aesthetic reasons, we want the fixed locus of $\psi^3$ to be empty, but we still
want to have the same distribution of the signs before the $x^i$. Therefore, we permute 
$\psi^2$ and $\psi^3$ as well as the second and third coordinate. We obtain:

\begin{equation}
\begin{array}{rcl}
\label{K3T3Case2}
\psi^1((x^1,x^2,x^3) +  \mathbb{Z}^3) & := & (x^1,-x^2,- x^3) +  \mathbb{Z}^3 \\
\psi^2((x^1,x^2,x^3) +  \mathbb{Z}^3) & := &  (-x^1,x^2,\tfrac{1}{2}-x^3) + \mathbb{Z}^3 \\
\psi^3((x^1,x^2,x^3) +  \mathbb{Z}^3) & := & (-x^1,-x^2,\tfrac{1}{2}+x^3) +  \mathbb{Z}^3 \\
\end{array}
\end{equation}

In the third case, we set both translation parts to zero and obtain:

\begin{equation}
\begin{array}{rcl}
\label{K3T3Case3}
\psi^1((x^1,x^2,x^3) + \mathbb{Z}^3) & := & (x^1,-x^2,- x^3) +  \mathbb{Z}^3 \\
\psi^2((x^1,x^2,x^3) +  \mathbb{Z}^3) & := & (-x^1,x^2,-x^3) +  \mathbb{Z}^3 \\
\psi^3((x^1,x^2,x^3) +  \mathbb{Z}^3) & := & (-x^1,-x^2,x^3) +  \mathbb{Z}^3 \\
\end{array}
\end{equation}

We choose $\rho^1,\rho^2:S\rightarrow S$ as one of the pairs that we have 
found in Section \ref{NonSymplecticZ22}. We can combine any of the 
three actions on $T^3$ with any non-symplectic $\mathbb{Z}^2_2$-action 
on $S$ and obtain a large number of quotients $(S\times T^3)/
\mathbb{Z}_2^2$.

\subsection{The first case} \label{The_first_case}
We investigate the three cases separately and start with the first
one where the action on $T^3$ is given by (\ref{K3T3Case1}).  
Our first step is to determine the fixed loci of all $\psi^i$ with $i=1,2,3$. 
Since $\psi^1$ maps $x^1$ to $x^1+\tfrac{1}{2}$, it has no fixed points. 
The same argument can be made for $\psi^3$. We denote the fixed locus 
of a map $f:X\rightarrow X$ by $\text{Fix}(f)$ and obtain:

\begin{equation}
\begin{array}{rcl}
\text{Fix}(\psi^1) & = & \emptyset \\
\text{Fix}(\psi^2) & = & \bigcup_{\epsilon_1,\epsilon_2\in \{0,1\}} 
\{(\tfrac{1}{2}\epsilon_1,x^2,\tfrac{1}{2}\epsilon_2) + \mathbb{Z}^3 | x^2\in\mathbb{R}\} \\
\text{Fix}(\psi^3) & = & \emptyset \\
\end{array}
\end{equation}

$\text{Fix}(\psi^2)$ consists of four circles and $\text{Fix}(\rho^2)$ 
is a disjoint union of complex curves. The connected components of 
$\text{Fix}(\rho^2\times\psi^2)$ are thus diffeomorphic two $S^1\times \Sigma$ where 
$\Sigma$ is a complex curve. $\rho^1\times \psi^1$ maps one of these connected 
components to another one and the singular locus of 
$(S\times T^3)/\mathbb{Z}^2_2$ therefore consists of two copies of 
$\text{Fix}(\rho^2)\times S^1$. Since the differential of a non-symplectic 
involution at a fixed point can be written as $\text{diag}(1,-1)\in
\mathbb{C}^{2\times 2}$, we can conclude that the fibers of the normal 
bundle of the singular locus are isomorphic to $\mathbb{C}^2/\{\pm \text{Id}\}$. 
The circle factor of the singular locus is parameterized by $x^2$ and $dx^2$ 
is thus a harmonic nowhere vanishing one-form on the singular locus. 
Therefore, all conditions of Theorem \ref{Thm-Joyce-Karigiannis} are satisfied 
and we obtain a smooth manifold $M$ with a torsion-free $G_2$-structure.  

We have to check if the fundamental group of $M$ is finite since this would 
guarantee that the holonomy of $M$ is the whole group $G_2$. We make use
of a theorem of Armstrong \cite{Armstrong} that is cited below.

\begin{Th}
\label{Pi1Quotient}
Let $G$ be a discontinuous group of homeomorphisms of a path connected, 
simply connected, locally compact metric space $X$, and let $H$ be the 
normal subgroup of $G$ generated by those elements which have fixed points. 
Then the fundamental group of the orbit space $X/G$ is isomorphic to the factor
group $G/H$.
\end{Th}

In the case where $\text{Fix}(\rho^2)=\emptyset$ the singular locus of 
$(S\times T^3)/\mathbb{Z}^2_2$ is empty and $M$ is covered by $S\times T^3$. 
Since $T^3=\mathbb{R}^3/\mathbb{Z}^3$, $(S\times T^3)/\mathbb{Z}^2_2$ can 
be written as a quotient of $S\times \mathbb{R}^3$ by a semidirect product 
$\mathbb{Z}^3 \rtimes \mathbb{Z}^2_2$. Therefore, we have $\pi_1(M) \cong 
\mathbb{Z}^3 \rtimes \mathbb{Z}^2_2$. Since this case is not particularly 
interesting, we assume from now on that $\text{Fix}(\rho^2)$ is not empty. 
A calculation shows that the elements of $\mathbb{Z}^3 \rtimes \mathbb{Z}^2_2$
whose action on $S\times\mathbb{R}^3$ has fixed points are precisely the unit 
element and the $\alpha_2\beta$, where $\beta$ is a translation along a vector 
$\{(x,0,z)^{\top} | x,z\in\mathbb{Z}\}$. These elements generate a group that is 
isomorphic to $\mathbb{Z}^2 \rtimes \mathbb{Z}_2$. Therefore, we have 
$\pi_1((S\times T^3)/\mathbb{Z}^2_2) \cong \mathbb{Z}\rtimes \mathbb{Z}_2$. 
Because of Corollary \ref{JoyceKarigiannisBetti} we also have $\pi_1(M) \cong 
\mathbb{Z}\rtimes \mathbb{Z}_2$ and the holonomy of $M$ is not the whole group 
$G_2$. 

In fact, the manifolds that we obtain in this case are an example
of \emph{barely $G_2$-manifolds} \cite{Grigorian,HarveyMoore}. 
A barely $G_2$-manifold is a quotient $(X\times S^1)/\mathbb{Z}_2$, 
where $X$ is a Calabi-Yau threefold and $\mathbb{Z}_2$ denotes a free 
group action that acts anti-holomorphically on $X$ and orientation-reversing
on the circle $S^1$. The holonomy of a barely $G_2$-manifold is 
$SU(3)\rtimes\mathbb{Z}_2$, which acts irreducibly on the tangent space. 
In our case, $X$ can be constructed as follows. The product of $S$ with 
the complex structure $I_2$ and the torus $T^2$ with coordinates 
$x^1$ and $x^3$ is a complex manifold. The projection of $\rho^2\times\psi^2$ to
$S\times T^2$ is a holomorphic involution that is non-symplectic on 
$(S,I_2)$. Since $\psi^2$ fixes $4$ points of $T^2$, the quotient 
$(S\times T^2)/\mathbb{Z}_2$ has $A_1$-singularities along $4$
copies of $\text{Fix}(\rho^2)$. These singularities can be blown
up. Since $S\times T^2$ has trivial canonical bundle, the blow-up is 
a Calabi-Yau manifold. This construction is known in the literature as the
Borcea-Voisin construction \cite{Borcea,Voisin}. The circle in the product 
$X\times S^1$ is parameterized by $x^2$ and the group $\mathbb{Z}_2$ 
by which we divide $X\times S^1$ is generated by the lift of $\rho^1\times\psi^1$ 
to the resolution $X\times S^1$ of $(S\times T^3)/\langle \rho^2
\times\psi^2\rangle$.   

We compute the Betti numbers of $M$. Since $\psi^1$ has no fixed points, 
$M':=(S\times T^3)/\langle \rho^1\times\psi^1 \rangle$ is a smooth manifold. Moreover, 
$(S\times T^3)/\mathbb{Z}^2_2 = M'/\langle \rho^2\times\psi^2 \rangle$ and we can 
determine the Betti numbers of $M$ by Corollary \ref{JoyceKarigiannisBetti}. 
In order to do this, we need the Betti numbers of the singular set $L$ of
$M'/\langle \rho^2\times\psi^2 \rangle$. We recall that in the generic case the 
fixed locus of $\rho^2$ consists of one complex curve 
of genus $\frac{1}{2}(22-r_2-a_2)$ and $\frac{1}{2}(r_2-a_2)$ rational curves. 
Therefore, we have

\begin{equation}
\begin{array}{rcl}
b^0(\text{Fix}(\rho^2)) & = & \frac{1}{2}(r_2-a_2) +1 \\
b^1(\text{Fix}(\rho^2)) & = & 22 - r_2 - a_2 \\ 
\end{array}
\end{equation}

These equations still hold in the exceptional case $(r_2,a_2,\delta_2) = 
(10,8,0)$ but not in the case where the fixed locus is empty that we 
have already excluded. Since $L$ consists of two copies of $S^1\times
\text{Fix}(\rho^2)$, we can conclude that 

\begin{equation}
\begin{array}{rcl}
b^0(L) = 2b^0(\text{Fix}(\rho^2)) & = & r_2-a_2 +2 \\
b^1(L) = 2b^0(\text{Fix}(\rho^2)) + 2b^1(\text{Fix}(\rho^2)_S) & = &  46 - r_2 - 3a_2 \\ 
\end{array}
\end{equation}

We also need the Betti numbers of $(S\times T^3)/\mathbb{Z}_2^2$. We compute
them by counting the harmonic forms on $S\times T^3$ that are invariant under the 
$\mathbb{Z}^2_2$-action. The harmonic $k$-forms on $T^3$ are precisely 
the $dx^{i_1\ldots i_k}$. By taking a look at (\ref{K3T3Case1}) we see that 
none of them, except $1$ and $dx^{123}$, is preserved by the 
$\mathbb{Z}^2_2$-action. Since $\rho^1$ and $\rho^2$ are 
involutions, the eigenvalues of their action on $H^2(S,\mathbb{C})$ 
are $\pm 1$. For $\epsilon_1,\epsilon_2 \in \{1,-1\}$ we define

\[
V_{\epsilon_1,\epsilon_2} :=
\{x\in H^2(S,\mathbb{C}) | (\rho^1)^{\ast}x= \epsilon_1x,\:
(\rho^2)^{\ast}x= \epsilon_2x  \}
\]

Since we always assume that $L_1\cap L_2 = \{0\}$, we have $V_{1,1}=\{0\}$
and  

\[
\dim{V_{1,-1}} + \dim{V_{-1,1}} + \dim{V_{-1,-1}} = \dim{H^2(S,\mathbb{C})} = 22
\] 

By a careful examination we see that there are no invariant harmonic $1$- or 
$2$-forms on $S\times T^3$ and that the invariant $3$-forms are precisely

\begin{itemize}
    \item $dx^{123}$, 
    \item $dx^1\wedge \beta_1$ with $[\beta_1]\in V_{1,-1}$ 
    \item $dx^2\wedge \beta_2$ with $[\beta_2]\in V_{-1,1}$ 
    \item $dx^3\wedge \beta_3$ with $[\beta_3]\in V_{-1,-1}$ 
\end{itemize} 

Thus, we have proven that the Betti numbers of the quotient 
$(S\times T^3)/\mathbb{Z}^2_2$ are $b^1=b^2=0$ and $b^3=23$. 
Therefore, it follows from Corollary \ref{JoyceKarigiannisBetti} that  

\begin{equation}
\begin{array}{rcl}
b^1(M) & = & 0  \\
b^2(M) & = & b^0(L) = r_2-a_2 +2 \\ 
b^3(M) & = & 23 + b^1(L) = 69 - r_2 - 3a_2 \\
\end{array}
\end{equation}

Finally, we create a table of the possible values of $(b^2,b^3)$. 
Since the minimal value of $r_1$ is $1$ and the minimal value of
$r_1+a_1$ is $2$, the restrictions $r_1 + r_2\leq 11$ and
$r_1+r_2+a_1+a_2 < 22$ from Theorem \ref{KoLe-Theorem} become  
$r_2\leq 10$ and $r_2+a_2 < 20$. We see that most pairs $(r_2,a_2)$ 
from Theorem \ref{Theorem_Z22} satisfy one of these restrictions except 
the $7$ pairs $(11,9)$, $(11,11)$, $(12,10)$, $(18,2)$, $(19,1)$, $(19,3)$
and $(20,2)$. Since there is a non-symplectic involution with invariants 
$(10,10,1)$, we have to include the pair $(r_2,a_2)=(10,10)$, too, 
although we have excluded the exceptional case $(10,10,0)$. We 
insert the allowed values of $(r_2,a_2)$ into the equations for $(b^2,b^3)$ and 
finally obtain the following theorem.

\begin{Th}
Let $\psi^i:T^3\rightarrow T^3$ with $i=1,2$ be the maps that we have defined 
in (\ref{K3T3Case1}). Moreover, let $\rho^i$ be the generators of 
one of the non-symplectic $\mathbb{Z}^2_2$-actions from Theorem 
\ref{Theorem_Z22} or \ref{KoLe-Theorem}. The maps $\rho^i\times \psi^i$ 
generate a group $\mathbb{Z}^2_2$ that acts on $S\times T^3$, where
$S$ is an appropriate K3 surface.  

The quotient $(S\times T^3)/\mathbb{Z}^2_2$ is a $G_2$-orbifold with
$A_1$-singularities. The singularities can be resolved and if the fixed locus
of $\rho^2$ is non-empty, we obtain a manifold $M$ with fundamental group 
$\mathbb{Z}\rtimes\mathbb{Z}_2$ and holonomy $SU(3)\rtimes\mathbb{Z}_2$. 
The Betti numbers of $M$ depend on the choice of the $\rho^i$ and can be 
found in the table below. 

\begin{center}
\begin{tabular}{l|l}
$b^2$ & $b^3$ \\

\hline

2 & 25, 29, 33, 37, 41, 45, 49, 53, 57, 61, 65 \\

\hline

4 & 27, 31, 35, 39, 43, 47, 51, 55, 59, 63, 67 \\

\hline

6 & 37, 41, 45, 49, 53, 57 \\ 

\hline

8 & 39, 43, 47, 51, 55 \\

\hline

10 & 41, 45, 49, 53, 57 \\ 

\hline

12 & 43, 47, 51, 55, 59 \\

\hline

14 & 45, 49 \\

\hline

16 & 47 \\ 

\hline

18 & 41, 45, 49 \\

\hline 

20 & 43, 47, 51 \\
\end{tabular}
\end{center}
\end{Th}

\subsection{The second case}  
We proceed to the case where the $\psi^i$ are defined by
(\ref{K3T3Case2}). The fixed loci of the $\rho^i\times\psi^i$ are 
given by 

\begin{equation}
\begin{array}{rcl}
\text{Fix}(\rho^1\times\psi^1) & := &  \text{Fix}(\rho^1) \times \bigcup_{\epsilon_1,\epsilon_2\in \{0,1\}} 
\{(x^1,\tfrac{1}{2}\epsilon_1,\tfrac{1}{2}\epsilon_2) + \mathbb{Z}^3 | x^1\in\mathbb{R}\} \\
\text{Fix}(\rho^2\times\psi^2) & := & \text{Fix}(\rho^2) \times \bigcup_{\epsilon_1,\epsilon_2\in \{0,1\}} 
\{(\tfrac{1}{2}\epsilon_1,x^2,\tfrac{1}{2}\epsilon_2 + \tfrac{1}{4}) + \mathbb{Z}^3 | 
x^2\in\mathbb{R}\} \\
\text{Fix}(\rho^3\times\psi^3) & := & \emptyset \\
\end{array}
\end{equation}

Since the last coordinate of $\text{Fix}(\psi^1)$ is a multiple of $\tfrac{1}{2}$
and the last coordinate of $\text{Fix}(\psi^2)$ is $\tfrac{1}{4}$ plus a multiple 
of $\tfrac{1}{2}$, the set of all points in $S\times T^3$ that are fixed by a 
$\rho^i\times\psi^i$ is the disjoint union

\begin{equation}
\label{FixedLoci}
\text{Fix}(\rho^1\times\psi^1)\:\dot{\cup}\: \text{Fix}(\rho^2\times\psi^2) 
\end{equation}

As in the previous subsection, the blow-up of $(S\times T^3)/\langle\rho^1\times\psi^1
\rangle$ along the singularities is a product $X\times S^1$ of a Calabi-Yau 
manifold of Borcea-Voisin type and a circle. $\rho^2\times\psi^2$ induces 
an involution $\imath$ of $X\times S^1$ that is anti-holomorphic on $X$. 
The difference to the previous case is that $\imath$ acts not freely on 
$X\times S^1$, at least if $\text{Fix}(\rho^2)\neq\emptyset$. Its fixed point 
set is diffeomorphic to $\text{Fix}(\rho^2\times\psi^2)/ \langle \rho^1\times\psi^1 \rangle$. 
$\psi^1$ maps each of the four circles in $\text{Fix}(\psi^2)$ to a distinct one. 
The singular locus of $(X\times S^1)/\langle\imath\rangle$ therefore consists 
of two copies of $\text{Fix}(\rho^2)\times S^1$. Each connected component of 
this set is the product of a complex curve and $S^1$. For the
same reasons as in the first case we can apply Theorem 
\ref{Thm-Joyce-Karigiannis} and construct a smooth manifold $M$ with a 
torsion-free $G_2$-structure out of $(S\times T^3)/\mathbb{Z}^2_2$.

The fundamental group can be computed by the same methods as in
the previous case. In the case where the fixed loci of $\rho^1$ and $\rho^2$
are non-empty, the elements of $\mathbb{Z}^3 \rtimes 
\mathbb{Z}^2_2$ that fix a vector in $\mathbb{R}^3$ are the unit element,
the $\psi^1\circ \beta$, where $\beta(x)=x+v$ with $v\in\{(0,y,z)^{\top} | 
y,z\in\mathbb{Z}\}$ and the $\psi^2\circ \beta$ with $v\in\{(x,0,z)^{\top} 
| x,z\in\mathbb{Z}\}$. These elements generate all of $\mathbb{Z}^3 
\rtimes \mathbb{Z}^2_2$ and therefore we have 

\begin{equation*}
\pi_1(M) \cong \pi_1((S\times T^3)/\mathbb{Z}^2_2) 
\cong (\mathbb{Z}^3 \rtimes \mathbb{Z}^2_2)/(\mathbb{Z}^3
\rtimes \mathbb{Z}^2_2) \cong \{1\}
\end{equation*}

If either $\text{Fix}(\rho^1)$ or $\text{Fix}(\rho^2)$ is empty, we obtain
as in Subsection \ref{The_first_case} $\pi_1(M)\cong \mathbb{Z}\rtimes
\mathbb{Z}_2$. $M$ can be written as $(X\times S^1)/\mathbb{Z}_2$, where the 
$\mathbb{Z}_2$-action is free but different from the one in Subsection
\ref{The_first_case}. Since the Betti numbers of $(S\times T^3)/
\mathbb{Z}^2_2$ and the shape of the singular locus are the same 
as in the previous subsection, we obtain barely $G_2$-manifolds 
with the same Betti numbers as before. 

We focus on the case where $\text{Fix}(\rho^1)$ and $\text{Fix}(\rho^2)$
are non-empty. Since $M$ is simply connected, it has indeed holonomy $G_2$. 
The Betti numbers of $M$ can be computed with help of the equation for
the Betti numbers of a blow-up and Corollary \ref{JoyceKarigiannisBetti}.
We see that $\text{Fix}(\rho^1\times\psi^1)$ and $\text{Fix}(\rho^2\times\psi^2)$
contribute the same to the cohomology of $M$. We denote the singular set of 
$(S\times T^3)/\mathbb{Z}^2_2$ by $L$ and obtain

\begin{equation*}
\begin{array}{rcl}
b^0(L) & = & b^0(\text{Fix}(\rho^1\times\psi^1)/\langle \rho^2\times\psi^2\rangle) +
b^0(\text{Fix}(\rho^2\times\psi^2)/\langle \rho^1\times\psi^1\rangle) \\[1mm]
& = & 2b^0(\text{Fix}(\rho^1)) + 2b^0(\text{Fix}(\rho^2))\\[1mm]  
& = & 4 + r_1 + r_2 - a_1 -a_2 \\[1mm]
b^1(L) & = & b^1(\text{Fix}(\rho^1\times\psi^1)/\langle \rho^2\times\psi^2\rangle) +
b^1(\text{Fix}(\rho^2\times\psi^2)/\langle \rho^1\times\psi^1\rangle) \\[1mm]
& = & 2b^0(\text{Fix}(\rho^1)) + 2b^0(\text{Fix}(\rho^2)) +  
2b^1(\text{Fix}(\rho^1)) + 2b^1(\text{Fix}(\rho^2))\\[1mm]
& = & 92 - r_1 - r_2 - 3a_1 - 3a_2 \\ 
\end{array}
\end{equation*}

The quotient $(S\times T^3)/\mathbb{Z}^2_2$ has the same Betti numbers
as in the previous case and we obtain:   

\begin{equation}
\label{BettiCase2}
\begin{array}{rcl}
b^1(M) & = & 0  \\
b^2(M) & = & b^0(L) = 4 + r_1 + r_2 - a_1 -a_2  \\ 
b^3(M) & = & 23 + b^1(L) = 115 - r_1 - r_2 - 3a_1 - 3a_2 \\
\end{array}
\end{equation}

\begin{Ex}
\label{JoyceKarigiannisExample1}
A special case of the above construction can be found in Example 7.2. in
\cite{JoyceKarigiannis}. The authors construct a K3 surface $X$ as a double 
cover of $\mathbb{CP}^2$ that is branched along a sextic curve $C$. 
The map $\alpha$ that permutes both branches is a non-symplectic involution 
with fixed locus $C$. Complex conjugation on $\mathbb{CP}^2$ induces a 
second involution $\beta$ of $X$ that has fixed locus $S^2$. $\alpha$ and 
$\beta$ commute and generate a non-symplectic $\mathbb{Z}^2_2$-action.
The fixed locus of $\alpha\beta$ is empty.
The authors extend the $\mathbb{Z}^2_2$-action on $X$ to $T^3\times X$
in the same way as we did in (\ref{K3T3Case2}). By resolving the 
singularities of $(T^3\times X)/\langle \alpha,\beta \rangle$, the authors 
obtain a $G_2$-manifold with Betti numbers $b^2 = 4$ and $b^3 = 67$. 

With help of Theorem \ref{FixedLocusTheorem} and the list of all triples 
$(r,a,\delta)$ we can deduce that that the invariants of $\alpha$ are 
$(1,1,1)$, the invariants of $\beta$ are $(11,11,1)$ and the invariants 
of $\alpha\beta$ are $(10,10,0)$. Moreover, the fixed lattices of $\alpha$ 
and $\beta$ are orthogonal to each other since the first one is spanned 
by the cohomology class of the curve $C$ and $\beta$ acts 
orientation-reversing on $C$. Therefore, all assumptions that we have made
in this subsection are satisfied and we obtain the same Betti numbers
as in \cite{JoyceKarigiannis}. 
\end{Ex}

\begin{Rem}
Joyce and Karigiannis \cite{JoyceKarigiannis} observe that 
the above example is in fact a twisted connected sum. In order to see 
this, they choose the metric on $T^3$ as $dx^1 + dx^2 + R^2 dx^3$, where 
$R$ is large. After that, they divide by $\mathbb{Z}^2_2$ and resolve the 
singularities. Finally, they cut the quotient $(T^3\times X)/\mathbb{Z}^2_2$ 
along the middle of the circle in the $x^3$-direction. In the limit 
$R\rightarrow\infty$, they obtain two asymptotically cylindrical parts. The 
first one contains the fixed locus of $\alpha$ and the second one contains 
the fixed locus of $\beta$. These two parts are in fact the same as we 
obtain by the Kovalev-Lee construction with the non-symplectic involutions 
$\alpha$ and $\beta$ as starting point. They are glued together such that 
the direct sum of their fixed lattices is primitively embedded into $L$. Of 
course, the above arguments can not only be made for the example from 
\cite{JoyceKarigiannis} but for any $\rho^i\times\psi^i$ such that $(\rho^1,
\rho^2)$ is a non-symplectic $\mathbb{Z}^2_2$-action and the $\psi^i$ 
are as in (\ref{K3T3Case2}). Therefore, all examples that we obtain in this 
subsection are in fact twisted connected sums. This explains why our 
formulas for $b^2$ and $b^3$ are precisely the same in Theorem 5.7.(c) 
in \cite{KoLe}. 
\end{Rem}

We compute the Betti numbers of our examples and check if we obtain 
$G_2$-manifolds with new values of $(b^2,b^3)$. We assume that the 
non-symplectic $\mathbb{Z}^2_2$-action on $S$ is generated by simple 
involutions since the $G_2$-manifolds that we obtain by Theorem 
\ref{KoLe-Theorem} and the $\psi^i$ from (\ref{K3T3Case2}) can already
be found in \cite{KoLe}. The Betti numbers depend only on the sums
$r_1 + r_2$ and $a_1 + a_2$. We take a look at the tables from Theorem 
\ref{Theorem_Z22} and see that $(r'_1+r'_2,a'_1+a'_2)$ is one of the 
following

\begin{align*}
& (2,2)\quad (3,1)\quad (3,3)\quad (4,0) \\
& (4,2)\quad (4,4)\quad (5,1)\quad (5,3) \\
& (5,5)\quad (6,4)\quad (6,6)\quad (7,5) \\
\end{align*}

Analogously, the pair $(r''_1+r''_2,a''_1+a''_2)$ is an element of the
following list:

\[
(0,0)\quad (8,0)\quad (8,8)\quad (16,0)\quad  (16,16)
\]

We obtain $60$ pairs $(r_1+r_2,a_1+a_2)$. $31$ of them satisfy
the condition $r_1+r_2\leq 11$ or $r_1+r_2+a_1+a_2 < 22$. Since these
cases are already studied in \cite{KoLe}, we do not consider them further. 
The Betti numbers of the remaining $29$ pairs can be calculated by the
equations (\ref{BettiCase2}). All in all, we have proven the following theorem. 

\begin{Th}
Let $\psi^i:T^3\rightarrow T^3$ with $i=1,2$ be the maps that we have defined 
in (\ref{K3T3Case2}). Moreover, let $\rho^i$ be the generators of 
one of the non-symplectic $\mathbb{Z}^2_2$-actions from Theorem 
\ref{Theorem_Z22} such that the fixed loci of the $\rho^i$ are non-empty.
The maps $\rho^i\times \psi^i$ generate a group $\mathbb{Z}^2_2$ that 
acts on $S\times T^3$, where $S$ is an appropriate K3 surface.  

The quotient $(S\times T^3)/\mathbb{Z}^2_2$ is a $G_2$-orbifold with
$A_1$-singularities. The singularities can be resolved and we obtain
a simply connected manifold $M$ with holonomy $G_2$. $M$ is 
diffeomorphic to the twisted connected sum that can be obtained from the
construction of Kovalev and Lee \cite{KoLe} with the non-symplectic involutions 
$\rho^1$ and $\rho^2$ as input. The above construction yields $G_2$-manifolds 
with $29$ distinct pairs of Betti numbers $(b^2,b^3)$ that are not included in 
Theorem 5.7.(c) in \cite{KoLe}. Their Betti numbers together with the invariants 
of the pairs of non-symplectic involutions can be found in the table below. 

\begin{center} 
\begin{tabular}{ccc}
\begin{tabular}{r|r|r|r}
$b^2$ & $b^3$ & $r_1+r_2$ & $a_1+a_2$ \\
\hline
4 & 27 & 22 & 22 \\ 
4 & 31 & 21 & 21 \\ 
4 & 35 & 20 & 20 \\ 
4 & 39 & 19 & 19 \\ 
4 & 43 & 18 & 18 \\ 
4 & 59 & 14 & 14 \\
4 & 63 & 13 & 13 \\ 
4 & 67 & 12 & 12 \\ 
6 & 29 & 23 & 21 \\ 
6 & 33 & 22 & 20 \\ 
6 & 37 & 21 & 19 \\ 
6 & 41 & 20 & 18 \\ 
6 & 45 & 19 & 17 \\ 
6 & 61 & 15 & 13 \\ 
6 & 65 & 14 & 12 \\ 
\end{tabular}
& \quad &
\begin{tabular}{r|r|r|r}
$b^2$ & $b^3$ & $r_1+r_2$ & $a_1+a_2$ \\
\hline
6 & 69 & 13 & 11 \\ 
6 & 73 & 12 & 10 \\ 
8 & 43 & 21 & 17 \\ 
8 & 47 & 20 & 16 \\ 
8 & 75 & 13 & 9 \\ 
20 & 75 & 22 & 6 \\ 
20 & 79 & 21 & 5 \\ 
20 & 83 & 20 & 4 \\ 
20 & 87 & 19 & 3 \\ 
22 & 77 & 23 & 5 \\ 
22 & 81 & 22 & 4 \\ 
22 & 85 & 21 & 3 \\ 
22 & 89 & 20 & 2 \\ 
24 & 91 & 21 & 1 \\
& & & \\ 
\end{tabular}
\end{tabular} 
\end{center}
\end{Th}

\begin{Rem}
The above theorem yields probably more than one diffeomorphism type of
$G_2$-manifolds for some of the pairs $(b^2,b^3)$. This may happen for
pairs of non-symplectic involutions with the same values of $r_1+r_2$
and $a_1+a_2$ but different invariants $(r_i,a_i,\delta_i)$. 
\end{Rem}

In the literature, there exist several examples of $G_2$-manifolds whose
Betti numbers are the same as one of the pairs $(b^2,b^3)$ in the above 
table. 

\begin{itemize}
    \item $(4,59)$, $(4,63)$, $(4,67)$, $(6,61)$, $(6,65)$, $(6,69)$, $(6,73)$,
    $(8,75)$, $(20,75)$, $(20,79)$, $(20,83)$, $(20,87)$ can be found in 
    Table 1 in \cite{KoLe}. Those examples are twisted connected sums, too,
    but at least one of the parts that are glued together is not constructed with 
    help of a non-symplectic involution. 
 
     \item $(4,35)$, $(6,41)$, $(8,47)$, $(22,89)$  can be found in Section 6.1
    in \cite{KoLe}. The corresponding $G_2$-manifolds are constructed from 
    pairs of non-symplectic involutions as in this section.  
 
    \item $(4, 43)$ can be found in Section 6.3 in \cite{KoLe} and is a twisted
    connected sum of another type. 

    \item $(4,27)$  and $(6,33)$ can be found in Section 12.4 in Joyce \cite{Joyce}. 

    \item $(4,39)$, $(6,45)$ can be found in Section 12.5 in \cite{Joyce}. 

    \item $(8,43)$ can be found in Section 12.7 in \cite{Joyce}. 
\end{itemize}

The remaining pairs

\[
(4,31)\quad (6,29)\quad (6,37)\quad (22,77)\quad (22,81)\quad (22,85)\quad (24, 91)   
\]

belong to new $G_2$-manifolds, or at least they do not appear in \cite{Joyce1,Joyce,KoLe}.

\subsection{The third case} \label{The_third_case}
Finally, we investigate the case where the maps $\psi^i$ are defined by
(\ref{K3T3Case3}). The fixed loci of the $\rho^i\times\psi^i$ are:  

\begin{equation}
\begin{array}{rcl}
\text{Fix}(\rho^1\times\psi^1) & = & \text{Fix}(\rho^1) \times 
\bigcup_{\epsilon_1,\epsilon_2\in \{0,1\}} 
\{(x^1,\tfrac{1}{2}\epsilon_1,\tfrac{1}{2}\epsilon_2) + \mathbb{Z}^3 | x^1\in\mathbb{R}\} \\
\text{Fix}(\rho^2\times\psi^2) & = & \text{Fix}(\rho^2) \times
\bigcup_{\epsilon_1,\epsilon_2\in \{0,1\}} 
\{(\tfrac{1}{2}\epsilon_1,x^2,\tfrac{1}{2}\epsilon_2) + \mathbb{Z}^3 | x^2\in\mathbb{R}\} \\
\text{Fix}(\rho^3\times\psi^3) & = & \text{Fix}(\rho^3) \times
\bigcup_{\epsilon_1,\epsilon_2\in \{0,1\}} 
\{(\tfrac{1}{2}\epsilon_1,\tfrac{1}{2}\epsilon_2,x^3) + \mathbb{Z}^3 | x^3\in\mathbb{R}\} \\
\end{array}
\end{equation}

It may happen that $\rho^1$ and $\rho^2$ have a common fixed point $p\in S$.
In this situation, $p$ is a fixed point of $\rho^3$, too. The points
$(\tfrac{1}{2}\epsilon_1,\tfrac{1}{2}\epsilon_2,\tfrac{1}{2}\epsilon_3) + 
\mathbb{Z}^3\in T^3$ with $\epsilon_i\in \{0,1\}$ are fixed by all $\psi^i$. 
Therefore, $S\times T^3$ has points that are fixed by the whole group
$\mathbb{Z}^2_2$. This means that the orbifold $(S\times T^3)/\mathbb{Z}^2_2$ 
has more complicated singularities of type $\mathbb{R}^7/\mathbb{Z}^2_2$ 
and we cannot apply Theorem \ref{Thm-Joyce-Karigiannis} directly. 

In principle, it is possible to resolve the singularities in two steps. The first
step would be to blow up the singularities of $(S\times T^3)/\mathbb{Z}_2$,
where $\mathbb{Z}_2$ is generated by $\rho^1\times\psi^1$. We obtain
a smooth manifold $N$ and $\rho^2\times \psi^2$ can be lifted to an 
involution $\imath$ of $N$. It can be shown that the singular set of 
$N/\langle\imath\rangle$ consists of $A_1$-singularities along disjoint
three-dimensional submanifolds. Since blowing up means in this context that
the metric is perturbed, it is not a priori clear that there exist harmonic 1-forms 
without zeroes on the three-dimensional submanifolds. Although this is probably 
true, we would need an analytical argument to show that this is indeed the case. 
Since this investigation is beyond the scope of this paper, we assume from 
now on that the fixed locus of $\rho^1$ is empty. This ensures that $\rho^2$ 
and $\rho^3$ have no common fixed point. 

In order to construct explicit examples, we have to know that the fixed locus of 
$\rho^1$ is empty and we need some information on the fixed loci of $\rho^2$
and $\rho^3$, too. Theorem \ref{KoLe-Theorem} contains only the invariants of $\rho^1$ 
and $\rho^2$. Since it is not straightforward to determine the invariants of 
$\rho^1\rho^2$, we restrict ourselves in this subsection to pairs $(\rho^1,\rho^2)$ 
from Theorem \ref{Theorem_Z22}. 

We show how to resolve the singularities in our situation. The fixed locus of
$\rho^2\times \psi^2$ consists of four copies of $\text{Fix}(\rho^2)\times S^1$, 
where $S^1$ has the coordinate $x^2$. Analogously, the fixed locus of $\rho^3
\times \psi^3$ consists of four copies of $\text{Fix}(\rho^3)\times S^1$, where 
$S^1$ has the coordinate $x^3$. Let $\Sigma \subseteq S$ be one of the 
connected components of $\text{Fix}(\rho^2)$. $\rho^1$ either maps
$\Sigma$ to itself such that $\rho^1:\Sigma\to\Sigma$ has no fixed points or it 
maps $\Sigma$ to another connected component $\Sigma'$. $\rho^1$ does not 
preserve the orientation of $\Sigma$ since it is anti-holomorphic with respect to 
$I_2$. The map $\psi^1$ reverses the orientation of the circle factor. 
$\text{Fix}(\rho^2\times\psi^2)/\langle \rho^1\times \psi^1 \rangle$ thus is a 
quotient of $\text{Fix}(\rho^2\times\psi^2)$ by an involution that preserves the 
orientation. The same argument can be made for the action of $\rho^1\times
\psi^1$ on $\text{Fix}(\rho^3\times \psi^3)$. All in all, the singular locus
of $(S\times T^3)/\mathbb{Z}^2_2$ is the disjoint union of four copies of 
$(\text{Fix}(\rho^2)\times S^1)/\mathbb{Z}_2$ and four copies of $(\text{Fix}(\rho^3)
\times S^1)/\mathbb{Z}_2$.

Unfortunately, we cannot apply Theorem \ref{Thm-Joyce-Karigiannis}
to resolve the singularities since the harmonic one-form $dx^2$ on 
$\text{Fix}(\rho^2) \times S^1$ is not projected to a well-defined 1-form 
on $(\text{Fix}(\rho^2)\times S^1)/\mathbb{Z}_2$. The reason for this is 
that $\mathbb{Z}_2$ acts as $x^2\mapsto -x^2$ on $S^1$. Nevertheless, 
$dx^2$ is projected to a $Z$-twisted 1-form and we can apply Corollary 
\ref{Thm-Joyce-Karigiannis-Twisted} to obtain a smooth $G_2$-manifold $M$. 
The $\mathbb{Z}_2$-principal bundle $Z$ can be defined simply as $\pi: 
\text{Fix}(\rho^2\times \psi^2) \rightarrow \text{Fix}(\rho^2\times \psi^2)/
\mathbb{Z}_2$. Again, the same argument can be made for $dx^3$. 

In order to obtain the Betti numbers of $M$  we need the $Z$-twisted Betti 
numbers of $L:=\text{Fix}(\rho^2\times\psi^2)/\langle \rho^1\times\psi^1
\rangle$. Let $p\in L$ be arbitrary and let $q\in \text{Fix}(\rho^2\times\psi^2)$ 
be a point with $\pi(q)=p$. The fiber $\pi^{-1}(p)$ of the $\mathbb{Z}_2$-principal 
bundle consists of $q$ and $f(q)$, where $f$ is the restriction of $\rho^1\times
\psi^1$ to $\text{Fix}(\rho^2\times\psi^2)$. Therefore, the fiber of the bundle 
$\bigwedge^k T^{\ast} L \otimes_{\mathbb{Z}_2} Z$ over $p$ consists 
of all equivalence classes 

\[
\{(\alpha_p,q), (-\alpha_p,f(q)) \}
\]

with $\alpha_p\in \bigwedge^k T^{\ast}_p L$. This means that the sections of 
$\bigwedge^k T^{\ast} L \otimes_{\mathbb{Z}_2} Z$ correspond to the $k$-forms 
on $L$ that are invariant under $\alpha \mapsto -f^{\ast} \alpha$. The $Z$-twisted 
Betti numbers count the linearly independent harmonic $k$-forms with this property. 
The harmonic $0$-forms on $L$ are precisely the functions that are constant on 
each connected component. For any pair of distinct connected components that are 
mapped to each other by $f$ we find one function with our invariance property. Since 
$\text{Fix}(\rho^2\times \psi^2)$ consists of four copies of $\text{Fix}(\rho^2)\times S^1$, 
$b^0(L,Z)$ is two times the number of connected components of $\text{Fix}(\rho^2)$ 
that are not mapped to itself by $\rho^1$.  

We compute $b^1(L,Z)$.  Let $\Sigma\subseteq \text{Fix}(\rho^2)$ be a
complex curve of genus $g$. A basis of the space of harmonic $1$-forms on 
$\Sigma\times S^1$ is given by $dx^2$, $g$ holomorphic forms on $\Sigma$ 
and $g$ anti-holomorphic forms on $\Sigma$. We assume that $f(\Sigma) =
\Sigma$. $dx^2$ corresponds to a $Z$-twisted $1$-form on $(\Sigma\times 
S^1)/\mathbb{Z}_2$. Since $f$ is anti-holomorphic on $\Sigma$, it maps the $g$ 
holomorphic forms to $g$ anti-holomorphic forms and vice versa. We conclude 
that the $(+1)-$ and $(-1)$-eigenspace of $-f^{\ast}$ acting on the harmonic 
$1$-forms are of the same dimension. All in all, it follows that $b^1((\Sigma\times S^1)/
\mathbb{Z}_2,Z)=g+1$. Next, we assume that $f$ maps $\Sigma$ to another 
complex curve $\Sigma'$. We can extend any harmonic $1$-form on
$\Sigma\times S^1$ to $(\Sigma\cup \Sigma')\times S^1$ such that it 
is invariant under $-f^{\ast}$. Therefore, we obtain $b^1(((\Sigma \cup \Sigma') 
\times S^1)/ \mathbb{Z}_2,Z)=2g+1$. The fixed locus of $\rho^3$ yields an 
analogous contribution to the Betti numbers of $M$. 

From now on, we restrict ourselves to triples $(\rho^1,\rho^2,\rho^3)$ with
the property that neither $\text{Fix}(\rho^2)$ nor $\text{Fix}(\rho^3)$ contains 
a pair of curves with the same genus. This ensures that $\rho^1$ maps
each connected component of  $\text{Fix}(\rho^2)$ or $\text{Fix}(\rho^3)$
to itself. Under this assumption, we have $b^0(L,Z)=0$ and 
$b^1((\Sigma\times S^1)/\mathbb{Z}_2,Z)=g+1$ for each complex curve
$\Sigma$ of genus $g$ in $\text{Fix}(\rho^2)$ or $\text{Fix}(\rho^3)$. 

If $\rho^2$ is not of one of the two exceptional types, $\text{Fix}(\rho^2)$ 
consists of one curve of genus $\frac{1}{2}(22 - r_2 - a_2)$ and 
$\frac{1}{2}(r_2 - a_2)$ rational curves. We have to exclude the cases 
$r_2-a_2\geq 4$ and $r_2+a_2=22$ but $r_2\neq a_2$ from our considerations. 
In both cases, we would have two or more rational curves in the fixed locus. 
We have to exclude the exceptional case $(r_2,a_2,\delta_2)=(10,8,0)$, too, 
since in that case the fixed locus consists of two elliptic curves. 

We are now able to determine the Betti numbers of $M$. The fixed locus 
of $\psi^2$ contributes four copies of $(\text{Fix}(\rho^2)\times S^1)/
\mathbb{Z}_2$ to the singular locus. Therefore, we obtain if 
$\text{Fix}(\rho^2)\neq\emptyset$:  

\begin{equation*}
\begin{aligned}
b^0(L,Z) & = 0 \\
b^1(L,Z) & = 4\cdot \frac{22 - r_2 - a_2}{2} + 4 + 4\cdot\frac{r_2 - a_2}{2}
= 48 - 4a_2 \\
\end{aligned}
\end{equation*}

Since $(S\times T^3)/\mathbb{Z}^2_2$ has as usual the Betti numbers
$b^1=b^2=0$ and $b^3=23$, the Betti numbers of our $G_2$-manifold 
$M$ are given by 

\begin{equation}
\begin{array}{rcl}
b^1(M) & = & 0  \\
b^2(M) & = & 0 \\ 
b^3(M) & = & 119 - 4(a_2+a_3)
\end{array}
\end{equation}

if the fixed loci of $\rho^2$ and $\rho^3$ are both non-empty. If in addition 
to the fixed locus of $\rho^1$ the fixed locus of $\rho^k$ with $k\in\{2,3\}$ is 
empty, we have $b^3(M)=71 - 4a_{5-k}$. It is impossible that all three fixed loci 
are empty. In that case, the rank of the fixed lattice of any $\rho^k:S \rightarrow S$ 
with $k\in\{1,2,3\}$ would be $10$. The direct sum of the fixed lattices 
would be of dimension $30$, which is more than the dimension of the K3 
lattice. 

Let $\mathbb{Z}^3 \rtimes \mathbb{Z}^2_2$ be the group that is generated
by translations with integer coefficients and the lifts of the $\rho^k\times \psi^k$ to 
isometries of $S\times \mathbb{R}^3$. The group that is generated by all elements 
of $\mathbb{Z}^3\rtimes\mathbb{Z}^2_2$ that have a fixed point is 
$\mathbb{Z}^3\rtimes\mathbb{Z}^2_2$ itself if $\rho^2$ and $\rho^3$ both have 
fixed points. Therefore, $M$ is simply connected and its holonomy is the whole 
group $G_2$. If only one fixed locus is non-empty, we obtain similarly as in the first 
case $\pi_1(M)=\mathbb{Z}\rtimes\mathbb{Z}_2$ and the holonomy is not the whole 
group $G_2$. Examples of this type are in fact barely $G_2$-manifolds 
$(X\times S^1)/\mathbb{Z}_2$, where $X$ is a Calabi-Yau manifold of Borcea-Voisin 
type. More precisely, $X$ is the blow-up of $(S\times T^3)/\langle \rho^k \times
\psi^k \rangle$, where $\rho^k$ is the non-symplectic involution with fixed points. 

We compute the possible values of $b^3$ in the case where the only
empty fixed locus is that of $\rho^1$. The invariants of $\rho^1$ are
$(10,10,0)$, which implies that $r'_1=2$, $a'_1=2$, $r''_1=8$ and
$r''_1=8$. We take a look at the table of Theorem \ref{Theorem_Z22}
and see that we necessarily have $a_2' + a_3'\in \{2,4,6\}$ and 
$a_2'' + a_3''=8$. Moreover, we can obtain all $3$ values of $a_2' + a_3'$ 
without choosing $\rho^2$ or $\rho^3$ as an involution with invariants 
$(10,10,0)$ or violating any of the further restrictions that we have imposed 
on $\rho^2$ and $\rho^3$. Explicit examples are given by $(r''_2,a''_2)=(8,8)$, 
$(r''_3,a''_3)=(0,0)$ and 

\[
((r'_2,a'_2),(r'_3,a'_3))\in \{((3,1),(1,1)),((3,3),(1,1)),((3,3),(3,3))\}
\]

Therefore, we
have $a_2+a_3\in\{10,12,14\}$ and obtain $G_2$-manifolds with

\[
b^3 \in \{63,71,79\}
\]

We turn to the case where the fixed loci of $\rho^1$ and $\rho^2$ are empty. 
In Section \ref{NonSymplecticZ22} we have seen that a simple non-symplectic
involution with invariants $(10,10,0)$ is up to conjugation the map whose 
restriction to $2(-E_8)$ maps $(x,y)$ to $(y,x)$ and whose restriction to $3H$ 
can be written as

\[
\begin{pmatrix}
& M_1 & \\
M_1 & & \\
& & M_2 \\
\end{pmatrix}
\] 

We assume without loss of generality that $\rho^1$ is precisely this map. $\rho^2$ 
has to commute with $\rho^1$ and shall have a fixed lattice that is isomorphic
to $H(2) \oplus (-E_8)(2)$ and is orthogonal to the fixed lattice of $\rho^1$. 
This is only possible if the restriction of $\rho^2$ to $2(-E_8)$ maps $(x,y)$ to
$(-y,-x)$. The restriction of $\rho^2$ to $3H$ has to be of type

\[
\begin{pmatrix}
& M_l & \\
M_l & & \\
& & M_2 \\
\end{pmatrix}
\]   

in order to satisfy all our conditions. The fixed lattice of $\rho^1$ is 

\[
\{ (v,v,0)^{\top} \in 3H | v\in H \}
\]

and the fixed lattice of $\rho^2$ is 

\[
\{ (v,M_l v,0)^{\top} \in 3H | v\in H \}
\]

The intersection of both lattices is trivial only if $l=2$. In this case,
the fixed lattice of $\rho^3 = \rho^1\rho^2$ is $H_3$ and thus we have 
$(r_3,a_3,\delta_3)=(2,0,0)$. It follows that we obtain only one barely 
$G_2$-manifold, which has $b^3 = 71$. All in all, we have proven the 
following theorem. 

\begin{Th}
Let $\psi^i:T^3\rightarrow T^3$ with $i=1,2$ be the maps that we have defined 
in (\ref{K3T3Case3}). Moreover, let $\rho^i$ be the generators of 
one of the non-symplectic $\mathbb{Z}^2_2$-actions from Theorem 
\ref{Theorem_Z22} such that the fixed locus of $\rho^1$ is empty and the
fixed loci of $\rho^2$ and $\rho^3=\rho^1\rho^2$ are non-empty.
Moreover, the fixed loci of $\rho^2$ and $\rho^3$ shall contain no
pair of complex curves with the same genus. The maps $\rho^i\times 
\psi^i$ generate a group $\mathbb{Z}^2_2$ that acts on $S\times T^3$, 
where $S$ is an appropriate K3 surface.  

The quotient $(S\times T^3)/\mathbb{Z}^2_2$ is a $G_2$-orbifold with
$A_1$-singularities. The singularities can be resolved and we obtain
a simply connected manifold $M$ with holonomy $G_2$. The second
Betti number of $M$ is always $0$ and the values of $b^3$ that can be obtained
by our construction are $63$, $71$, and $79$. 

If the fixed locus of $\rho^2$ is empty, $M$ is a barely $G_2$-manifold
of a unique diffeomorphism type. Its Betti numbers are $b^2=0$ and
$b^3=71$. 
\end{Th}

\begin{Rem}
\begin{enumerate}
    \item $G_2$-manifolds with $b^2=0$ and the three values of $b^3$ 
    from the above theorem can also be found in Table 1 in \cite{KoLe}
    and in the section "$G_2$-manifolds from pairs of smooth Fano 3-folds"
    in Corti et al. \cite{CortiEtAl2}.  

    \item A set of invariants that is sufficient to determine the diffeomorphism type
    of a 2-connected 7-dimensional spin manifold can be found in 
    \cite{CrowleyNordstroem}. Although it is outside of the scope of this paper,
    it would be interesting to examine if our $G_2$-manifolds are 2-connected,
    to determine their diffeomorphism type and to check if they coincide with
    the examples from \cite{CortiEtAl2}. 
\end{enumerate}    
\end{Rem}

\begin{Ex}
Let $X$ be the K3 surface with the two non-symplectic involutions $\alpha$
and $\beta$ from \cite{JoyceKarigiannis} that we have already considered in
Example \ref{JoyceKarigiannisExample1}. Let $\alpha$ and $\beta$ act on 
$X\times T^3$ as in (\ref{K3T3Case3}). The resolution of the quotient
$(X\times T^3)/\langle\alpha,\beta\rangle$ is studied in Example 7.3. in 
\cite{JoyceKarigiannis}. The authors obtain a smooth $G_2$-manifold with 
$b^2=0$ and $b^3=71$. We recall that the invariants of $\alpha$ are $(1,1,1)$,
the invariants of $\beta$ are $(11,11,1)$ and that the fixed locus of $\alpha\beta$ is
empty. By our formula we obtain the same value of $b^3$. 
\end{Ex}

\subsection{Quotients by a non-abelian group}
Finally, we study resolutions of quotients of $S\times T^3$ by a non-abelian group 
$\Gamma$. We define two maps $\psi^1,\psi^2: T^3\rightarrow T^3$ by 

\begin{equation}
\label{K3T3Case4}
\begin{array}{rcl}
\psi^1((x^1,x^2,x^3) + \mathbb{Z}^3) & := & (x^1 + \tfrac{1}{4},-x^2 + \tfrac{1}{4},-x^3) + \mathbb{Z}^3 \\
\psi^2((x^1,x^2,x^3) + \mathbb{Z}^3) & := & (-x^1,x^2,-x^3) + \mathbb{Z}^3 \\
\end{array}
\end{equation}

Since 

\[
(\psi^1)^2((x^1,x^2,x^3) + \mathbb{Z}^3) = (x^1 + \tfrac{1}{2}, x^2, x^3) + \mathbb{Z}^3\:,
\]

we have $(\psi^1)^4 = (\psi^2)^2 = \text{Id}_{T^3}$. Moreover, we have  

\[
\psi^1\psi^2((x^1,x^2,x^3) + \mathbb{Z}^3) = (-x^1 + \tfrac{1}{4},-x^2 + \tfrac{1}{4},x^3) 
+  \mathbb{Z}^3 
\]  

and thus $(\psi^1\psi^2)^2= \text{Id}_{T^3}$. All in all, $\psi^1$ and $\psi^2$ 
generate a group $\Gamma$ that is isomorphic to the dihedral group $D_4$. We extend
the $\psi^i$ to maps $\rho^i\times\psi^i:S\times T^3 \to S\times T^3$ by assuming that 

\begin{equation}
\begin{aligned}
& {\rho^1}^{\ast}\omega_1 = \omega_1\:, 
\qquad 
{\rho^1}^{\ast}\omega_2 = - \omega_2\:,
\qquad
{\rho^1}^{\ast}\omega_3 = -\omega_3 \\
& {\rho^2}^{\ast}\omega_1 = - \omega_1\:, 
\qquad 
{\rho^2}^{\ast}\omega_2 = \omega_2\:,
\qquad
{\rho^2}^{\ast}\omega_3 = - \omega_3 
\end{aligned}
\end{equation}

We determine the fixed loci of all elements of $\Gamma$. For reasons of brevity, we
denote these elements by  

\[
\gamma_{jk} := (\rho^1\times\psi^1)^j (\rho^2\times \psi^2)^k 
\]

with $j\in\{0,1,2,3\}$ and $k\in\{0,1\}$. $(\psi^1)^j$ with $j\in\{1,2,3\}$ acts on the 
$x^1$-coordinate as a translation. Therefore, the fixed point set of $\gamma_{j0}$ 
is all of $S\times T^3$ if $j=0$ and empty if $j\neq 0$. By a short calculation we 
see that the fixed point sets of the remaining group elements are given by:

\begin{equation}
\begin{array}{rcl}
\text{Fix}(\gamma_{01}) & = & \bigcup_{\epsilon_1,\epsilon_2\in \{0,1\}} 
\{(\tfrac{1}{2}\epsilon_1,x^2,\tfrac{1}{2}\epsilon_2) + \mathbb{Z}^3 | x^2\in\mathbb{R}\}
\times \text{Fix}(\rho^2) \\
\text{Fix}(\gamma_{11}) & = & \bigcup_{\epsilon_1,\epsilon_2\in \{0,1\}} 
\{(\tfrac{1}{8} + \tfrac{1}{2}\epsilon_1,\tfrac{1}{8} + \tfrac{1}{2}\epsilon_2,x^3) + 
\mathbb{Z}^3 | x^3\in\mathbb{R}\} \times \text{Fix}(\rho^1\rho^2) \\
\text{Fix}(\gamma_{21}) & = & \bigcup_{\epsilon_1,\epsilon_2\in \{0,1\}} 
\{(\tfrac{1}{4} + \tfrac{1}{2}\epsilon_1,x^2, \tfrac{1}{2}\epsilon_2) + \mathbb{Z}^3 | x^2\in\mathbb{R}\}
\times \text{Fix}(\rho^2) \\
\text{Fix}(\gamma_{31}) & = & \bigcup_{\epsilon_1,\epsilon_2\in \{0,1\}} 
\{(\tfrac{3}{8} + \tfrac{1}{2}\epsilon_1,\tfrac{1}{8} + \tfrac{1}{2}\epsilon_2,x^3) + 
\mathbb{Z}^3 | x^3\in\mathbb{R}\} \times \text{Fix}(\rho^1\rho^2) \\
\end{array}
\end{equation}

$\rho^1\times\psi^1$ maps $\text{Fix}(\gamma_{01})$ bijectively to $\text{Fix}(\gamma_{21})$
and $\text{Fix}(\gamma_{11})$ to $\text{Fix}(\gamma_{31})$. The singular locus of $(S\times 
T^3)/\Gamma$ therefore consists of one component that is diffeomorphic to 
$\text{Fix}(\gamma_{01})$ and another one that is diffeomorphic to $\text{Fix}(\gamma_{11})$.
More precisely, it is the union of $4$ copies of $S^1\times \text{Fix}(\rho^2)$ and $4$ 
copies of $S^1\times \text{Fix}(\rho^1\rho^2)$. By taking a look at the first coordinate of $T^3$ 
we see that this union is disjoint. Since $\gamma_{01}$ and $\gamma_{11}$ are both of order 
$2$, all singularities of $(S\times T^3)/\Gamma$ are $A_1$-singularities. Therefore, we can 
apply Theorem \ref{Thm-Joyce-Karigiannis} to obtain a smooth $G_2$-manifold $M$. 

We determine the fundamental group of $M$. As before, we have a natural action of 
$\mathbb{Z}^3\rtimes\Gamma$ on $S\times \mathbb{R}^3$. We consider the following 
maps that are induced by this group action:

\begin{equation}
\begin{array}{rcl}
(p,(x^1,x^2,x^3)) & \mapsto & (\rho^2(p),(-x^1,x^2,-x^3) + (\lambda_1,0,\lambda_2)) \\
(p,(x^1,x^2,x^3)) & \mapsto & (\rho^1\rho^2(p),(-x^1+\frac{1}{4},-x^2+\frac{1}{4},x^3) + (\lambda_1,\lambda_2,0)) \\ 
\end{array}
\end{equation}

where $p\in S$, $(x^1,x^2,x^3)\in\mathbb{R}^3$ and $\lambda_1,\lambda_2\in\mathbb{Z}$
are arbitrary. If $\rho^2$ and $\rho^1\rho^2$ have at least one fixed point, the above maps 
have a fixed point, too, and we can conclude with help of Theorem \ref{Pi1Quotient} that 
$M$ is simply connected and thus has holonomy $G_2$. If $\text{Fix}(\rho^2)$ or 
$\text{Fix}(\rho^1\rho^2)$ is empty, the fundamental group is infinite and the holonomy
is smaller. We do not pursue this case further. 

Finally, we calculate the Betti numbers of $M$. By the same arguments as in the previous
cases we see that $b^1((S\times T^3)/\Gamma) = b^2((S\times T^3)/\Gamma) = 0$ and
$b^3((S\times T^3)/\Gamma) = 23$. We denote the invariants of $\rho^2$ by $(r_2,a_2,\delta_2)$ 
and the invariants of $\rho^3$ by $(r_3,a_3,\delta_3)$. The Betti numbers of the singular locus $L$ 
can be computed by the usual methods as follows:

\begin{equation*}
\begin{array}{rcl}
b^0(L) & = & 4\cdot\left(\frac{r_2 - a_2}{2} + 1 +  \frac{r_3 - a_3}{2} + 1\right) 
= 8 + 2(r_2 + r_3) - 2(a_2 + a_3) \\
b^1(L) & = & b^0(L) + 4\cdot\left((22 - r_2 - a_2) + (22 - r_3 - a_3) \right) \\
& = & 184 - 2(r_2 + r_3) - 6(a_2 + a_3) \\
\end{array}
\end{equation*}

The above equations yield the correct Betti numbers for 
$(r_k,a_k,\delta_k)=(10,8,0)$ with $k\in\{2,3\}$, too. We obtain the
following Betti numbers for $M$: 

\begin{equation}
\begin{array}{rcl}
b^2(M) & = & 8 + 2(r_2 + r_3) - 2(a_2 + a_3) \\
b^3(M) & = & 207 - 2(r_2 + r_3) - 6(a_2 + a_3) \\
\end{array}
\end{equation}

The Betti numbers of $M$ depend only on $r_2+r_3$ and $a_2+a_3$. The
set of all pairs $(r_2 + r_3,a_2 + a_3)$ from Theorem \ref{Theorem_Z22} and 
Theorem \ref{KoLe-Theorem} is the same as the set of all $(r_1 + r_2,a_1 + 
a_2)$ since $(\rho^2,\rho^3)$ generate the same non-symplectic 
$\mathbb{Z}^2_2$-action as $(\rho^1,\rho^2)$. We calculate the Betti numbers 
that we obtain from these pairs and finally have proven the following theorem. 

\begin{Th}
Let $\psi^i:T^3\rightarrow T^3$ with $i=1,2$ be the maps that we have defined 
in (\ref{K3T3Case4}). Moreover, let $\rho^i$ be the generators of 
one of the non-symplectic $\mathbb{Z}^2_2$-actions from Theorem 
\ref{Theorem_Z22} or \ref{KoLe-Theorem}. The maps $\rho^i\times \psi^i$ 
generate a group that is isomorphic to the dihedral group $D_4$ and acts on 
$S\times T^3$, where $S$ is an appropriate K3 surface.  

The quotient $(S\times T^3)/D_4$ is a $G_2$-orbifold with
$A_1$-singularities. The singularities can be resolved and we obtain
a manifold $M$ whose holonomy is contained in $G_2$. If the fixed loci
of $\rho^2$ and $\rho^1\rho^2$ are non-empty, $M$ is simply connected and
its holonomy is $G_2$. The Betti numbers of $M$ depend on the choice of 
the $\rho^i$ and can be found in the table below. 

\begin{center}
\begin{tabular}{l|l}
$b^2$ & $b^3$ \\

\hline

8 & 31, 39, 47, 55, 63, 95, 103, 111, 119, 127, 135, 143, 151, 159, 167,\\
& 175, 183, 191 \\

\hline

12 & 35, 43, 51, 59, 67, 99, 107, 115, 123, 131, 139, 147, 155, 163, 171,\\ 
& 179, 187, 195 \\

\hline

16 & 63, 71, 127, 135, 143, 151, 159, 167, 175, 183, 191, 199 \\

\hline

20 & 139, 147, 155, 163, 171, 179 \\

\hline

24 & 79, 87, 95, 103, 111, 143, 151, 159, 167, 175 \\

\hline

28 & 83, 91, 99, 107, 115, 147, 155, 163, 171, 179 \\

\hline

32 & 111, 119, 151, 159, 167, 175, 183 \\

\hline

36 & 155, 163 \\

\hline

40 & 127, 135, 143, 151, 159 \\

\hline 

44 & 131, 139, 147, 155, 163  \\

\hline

48 & 159, 167 \\ 
\end{tabular}
\end{center}
\end{Th}

The above theorem yields $95$ different pairs $(b^2,b^3)$. Some of them 
can be found in the literature:

\begin{itemize}
    \item $(8,47)$, $(8,63)$, $(8,95)$, $(8,103)$, $(8,111)$, $(8,119)$, $(8,151)$, $(8,159)$,
    $(12,59)$, $(12,67)$, $(12,99)$, $(12,107)$, $(12,115)$, $(12,123)$, $(12,155)$, $(12,163)$,
    $(16,71)$, $(20,155)$ and $(24,95)$ can be found in Kovalev, Lee \cite{KoLe}.  

    \item Of the remaining pairs, $(8,31)$, $(8,39)$, $(8,55)$, $(8,127)$, $(8,135)$, $(12,43)$, $(12,51)$,
    $(12,131)$ and $(16,135)$ can be found in Joyce \cite{Joyce}. 
\end{itemize}

The remaining $67$ pairs do not appear in \cite{Joyce1,Joyce,KoLe}. One pair of Betti 
numbers from the above theorem usually corresponds to several non-symplectic 
$\mathbb{Z}^2_2$-actions and thus probably to several diffeomorphism types $G_2$-manifolds. 
Although the pairs $(b^2,b^3)$ from the literature usually correspond to several $G_2$-manifolds, 
too, it is possible that our theorem yields further new $G_2$-manifolds.

\end{document}